\documentclass[12pt]{article}
\usepackage{amssymb}
\usepackage{amsmath}

\input{tcilatex}
\begin{document}

\begin{center}
\textbf{A generalization for a finite family of functions of the converse of
Browder's fixed point theorem}

\bigskip

Radu MICULESCU and Alexandru MIHAIL

\bigskip

Bucharest University, Faculty of Mathematics and Computer Science

Str. Academiei 14, 010014 Bucharest, Romania; miculesc@yahoo.com

\medskip

Bucharest University, Faculty of Mathematics and Computer Science

Str. Academiei 14, 010014 Bucharest, Romania; mihail\_alex@yahoo.com

\bigskip
\end{center}

\textbf{Abstract}. {\small Taking as model the attractor of an iterated
function system consisting of }$\varphi ${\small -contractions on a complete
and bounded metric space, we introduce the set-theoretic concept of family
of functions having attractor. We prove that, given such a family, there
exist a metric on the set on which the functions are defined and take values
and a comparison function }$\varphi $ {\small such that all the family's
functions are} $\varphi ${\small -contractions. In this way we obtain a
generalization for a finite family of functions of the converse of Browder's
fixed point theorem. As byproducts we get a particular case of Bessaga's
theorem concerning the converse of the contraction principle and a companion
of Wong's result which extends the above mentioned Bessaga's result for a
finite family of commuting functions with common fixed point.}

\bigskip



\bigskip

\textbf{1.} \textbf{INTRODUCTION}

\bigskip

The problem of the converse of Banach-Picard-Caccioppoli principle was
treated by several mathematicians each of them concentrating on different
assumptions. C. Bessaga (see [3], [10] and [13]) was the first one to treat
the problem by using only set-theoretic assumptions. J. S. W. Wong (see
[24]) extended Bessaga's result for a finite family of commuting functions
with common unique fixed point. Other results on this direction are due to
L. Jano\v{s} (see [12]), P.R. Meyers (see [18]) and S. Leader (see [15]).

The idea of replacing the contractivity condition imposed on the function $%
f:X\rightarrow X$ considered in the Banach-Picard-Caccioppoli principle by a
weaker one described by the inequality $d(f(x),f(y))\leq \varphi (d(x,y))$
for all $x,y\in X$, where $\varphi $ has certain properties defining the so
called comparison function, was treated, among others, by D.W. Boyd and J.S.
Wong (see [4]), F. Browder (see [5]), J. Matkowski (see [17]) and I. A. Rus
(see [21]). A function $f$ satisfying the previous inequality is called $%
\varphi $-contractions. From the point of view of the problem treated in
this paper a special place is played by Browder's result concerning $\varphi 
$-contractions (see Theorem 2.5). For more details about this result one can
consult [11].

Iterated function systems, introduced by J. Hutchinson (see [9]) and
popularized by M. Barnsley (see [1]), represent one of the most general way
to generate fractals. The large variety of their applications is the
background of the current effort to extend the classical Hutchinson's
theory. One line of research in this direction is to weaken the usual
contraction condition by considering iterated function systems consisting of 
$\varphi $-contractions. For results in this direction one can consult [7],
[8], [22] and [23].

By selecting some properties of the attractor of an iterated function system
consisting of $\varphi $-contractions on a complete and bounded metric
space, we introduced the set-theoretic concept of family of functions having
attractor (Definition 3.3). We prove that, given such a family, there exist
a complete and bounded metric on the set on which the functions are defined
and take values and a comparison function $\varphi $ such that all the
family's functions are $\varphi $-contractions (see Theorem 3.21). In this
way we obtain a generalization for a finite family of functions of the
converse of Browder's fixed point theorem.

If $\mathcal{F}=(f_{i})_{i\in I}$ is a family of functions having attractor $%
A$, where $f_{i}:X\rightarrow X$\ and $I$\ is finite, we obtain the result
tracking the following steps:

- the construction (based on the main result from [19]) of a metric $d$ on $%
A $ and a comparison function $\varphi $ such that $d(f_{i}(x),f_{i}(y))\leq
\varphi (d(x,y))$\textit{\ }for every $i\in I$ and every $x,y\in A$, i.e.%
\textit{\ }$f_{i}$'s are $\varphi $-contractions on the attractor with
respect to $d$ (Theorem 3.4)

- the construction of a semi-metric $d^{\mu }$\ on $X$, associated to $%
\mathcal{F}$ and to a sequence $\mu $, such that $d^{\mu
}(f_{i}(x),f_{i}(y))\leq d^{\mu }(x,y)$ for every $x,y\in X$, i.e. $f_{i}$'s
are nonexpansive on $X$ with respect to $d^{\mu }$ (Proposition 3.8)

- the construction of a complete and bounded metric $d$ on $X$ (Proposition
3.16)

- the construction of a comparison function $\varphi $\ such that $%
d(f_{i}(x),f_{i}(y))\leq \varphi (d(x,y))$\ for every $i\in I$ and every $%
x,y\in X$, i.e. $f_{i}$'s are $\varphi $-contractions with respect to $d$
(Lemma 3.20).

Finally we present a result which removes the boundedness condition on the
metric $d$, we point out that one can obtain from our result a particular
case of Bessaga's theorem concerning the converse of the contraction
principle (see Theorem 5 from [10]) and we present a companion of Wong's
result which extends the above mentioned Bessaga's result for a finite
family of commuting functions with common fixed point (see [24]).

\bigskip

\textbf{2.} \textbf{PRELIMINARIES}

\bigskip

For a function $f:X\rightarrow X$ and $n\in \mathbb{N}$, by $f^{[n]}$\ we
mean the composition of $f$\ by itself $n$\ times.

\bigskip

\textbf{Definition 2.1} (comparison function)\textbf{.} \textit{A function }$%
\varphi :[0,\infty )\rightarrow \lbrack 0,\infty )$\textit{\ is called a
comparison function if it has the following three properties:}

\textit{i) }$\varphi $\textit{\ is increasing;}

\textit{ii) }$\varphi (t)<t$\textit{\ for every }$t>0$\textit{;}

\textit{iii) }$\varphi $\textit{\ is right-continuous.}

\bigskip

\textbf{Remark 2.2.}

i) Any function $\varphi :[0,\infty )\rightarrow \lbrack 0,\infty )$\
satisfying ii) and iii) from the above definition has the following
property: $\underset{n\rightarrow \infty }{\lim }\varphi ^{\lbrack n]}(t)=0$
for every $t>0$ (see Remark 1 from [16]).

ii) $\varphi (0)=0$ for every comparison function.

\bigskip

\textbf{Definition 2.3 }($\varphi $-contraction)\textbf{.} \textit{Let }$%
(X,d)$\textit{\ be a metric space\ and a function }$\varphi :[0,\infty
)\rightarrow \lbrack 0,\infty )$\textit{. A function }$f:X\rightarrow X$%
\textit{\ is called a }$\varphi $\textit{-contraction if} $d(f(x),f(y))\leq
\varphi (d(x,y))$ \textit{for all }$x,y\in X$\textit{.}

\bigskip

\textbf{Remark 2.4.} Every $\varphi $-contraction is Lipschitz, so it is
continuous.

\bigskip

The next result is known as Browder's Theorem.

\bigskip

\textbf{Theorem 2.5 }(see Theorem 1 from [5], Theorem 1 from [11] or Example
2.9., 1) from [2])\textbf{.} \textit{Let }$(X,d)$\textit{\ be a complete and
bounded metric space\ and }$\varphi :[0,\infty )\rightarrow \lbrack 0,\infty
)$ \textit{a comparison function. Then every }$\varphi $\textit{-contraction 
}$f:X\rightarrow X$\textit{\ has a unique fixed point }$x_{0}$\textit{\ and }%
$\underset{n\rightarrow \infty }{\lim }f^{[n]}(x)=x_{0}$ \textit{for every }$%
x\in X$.

\bigskip

Given a metric space $(X,d)$ and a subset $Y$ of $X$, by $d(Y)$ we denote
the diameter of $Y$ and by $\mathcal{K}(X)$ we denote the family of
non-empty compact subsets of $X$.

\bigskip

For a nonempty set $I$, by $\Lambda (I)$ we mean the set $I^{\mathbb{N}%
^{\ast }}$ and by $\Lambda _{n}(I)$ we mean the set $I^{\{1,2,...,n\}}$. So,
the elements of $\Lambda (I)$ are written as infinite words $\alpha =\alpha
_{1}\alpha _{2}...\alpha _{m}\alpha _{m+1}...$ and the elements of $\Lambda
_{n}(I)$ are written as finite words $\alpha =\alpha _{1}\alpha
_{2}...\alpha _{n}$ ($n$, which is the length of $\omega $, is denoted by $%
\left\vert \omega \right\vert $).

By $\Lambda ^{\ast }(I)$ we denote the set of all finite words, i.e. $%
\Lambda ^{\ast }(I)\overset{def}{=}\underset{n\in \mathbb{N}^{\ast }}{\cup }%
\Lambda _{n}(I)\cup \{\lambda \}$, where $\lambda $ is the empty word.

For $\alpha =\alpha _{1}\alpha _{2}...\alpha _{m}\alpha _{m+1}...\in \Lambda
(I)$ and $n\in \mathbb{N}$, we shall use the following notation: $[\alpha
]_{n}\overset{not}{=}\alpha _{1}\alpha _{2}...\alpha _{n}$ if $n\geq 1$ and $%
\lambda $ if $n=0$.

For two words $\alpha \in \Lambda _{n}(B)\,$and $\beta \in \Lambda _{m}(B)$
or $\beta \in \Lambda (B)$, by $\alpha \beta $ we mean the concatenation of
the words $\alpha $ and $\beta $, i.e.$\ \alpha \beta =\alpha _{1}\alpha
_{2}...\alpha _{n}\beta _{1}\beta _{2}...\beta _{m}$ and respectively $%
\alpha \beta =\alpha _{1}\alpha _{2}...\alpha _{n}\beta _{1}\beta
_{2}...\beta _{m}\beta _{m+1}...$.

On $\Lambda (I)$ we consider the metric given by $d_{\Lambda }(\alpha ,\beta
)=\overset{\infty }{\underset{k=1}{\sum }}\frac{1-\delta _{\alpha
_{k}}^{\beta _{k}}}{3^{k}}$, where $\delta _{x}^{y}=\{%
\begin{array}{c}
1\text{, if }x=y \\ 
0\text{, if }x\neq y%
\end{array}%
$.

\bigskip

\textbf{Remark 2.6.} The function $\tau _{i}:\Lambda (I)\rightarrow \Lambda
(I)$, given by $\tau _{i}(\alpha )=i\alpha $ for every $\alpha \in \Lambda
(I)$, is continuous.

\bigskip

\textbf{Remark 2.7.}

i) The convergence in the compact metric space $(\Lambda (I),d_{\Lambda })$\
is the convergence on components.

ii) If $I$ is finite, then $(\Lambda (I),d_{\Lambda })$\ is compact.

\bigskip

Given the functions $f_{i}:X\rightarrow X$, where $X$ is a given set and $%
i\in I$, we shall use the following notations:

i) $f_{\lambda }=Id_{X}$;

ii) $f_{\alpha _{1}\alpha _{2}...\alpha _{m}}\overset{not}{=}f_{\alpha
_{1}}\circ f_{\alpha _{2}}\circ ...\circ f_{\alpha _{m}}$ for every $\alpha
_{1},\alpha _{2},...,\alpha _{m}\in I$;

iii) $Y_{\alpha }\overset{not}{=}f_{\alpha }(Y)$ for every $\alpha \in
\Lambda ^{\ast }(I)$ and every $Y\subseteq X$.

\bigskip

\textbf{Definition 2.8} (topological self-similar set, topological
self-similar system)\textbf{.} \textit{A compact Hausdorff topological space 
}$K$\textit{\ is called a topological self-similar set if there exist
continuous functions }$f_{1}$\textit{, }$f_{2}$\textit{, ..., }$%
f_{N}:K\rightarrow K$\textit{, where} $N\in \mathbb{N}^{\ast }$\textit{,} 
\textit{and a continuous surjection }$\pi :\Lambda
(\{1,2,...,N\})\rightarrow K$\textit{\ such that the diagram }%
\begin{equation*}
\begin{array}{ccc}
\text{ \ \ }\Lambda (\{1,2,...,N\}) & \overset{\tau _{i}}{\rightarrow } & 
\Lambda (\{1,2,...,N\}) \\ 
\pi \downarrow &  & \text{ \ }\downarrow \pi \\ 
\text{ \ \ }K & \underset{f_{i}}{\rightarrow } & K%
\end{array}%
\end{equation*}%
\textit{commutes for all }$i\in \{1,2,...,N\}$\textit{.}

\textit{We say that }$(K,(f_{i})_{i\in \{1,2,...,N\}})$\textit{, a
topological self-similar set together with the set of continuous maps as
above, is a topological self-similar system.}

\bigskip

The above definition is Definition 0.3 from [14].

\bigskip

\textbf{Theorem 2.9} (see Theorem 3.1 from [19])\textbf{.} \textit{For every
topological self-similar system }$(K,(f_{i})_{i\in \{1,2,...,N\}})$\textit{\
there exist a metric }$\delta $\textit{\ on }$K$\textit{\ which is
compatible with the original topology and a comparison function }$\varphi
:[0,\infty )\rightarrow \lbrack 0,\infty )$\textit{\ such that} $\delta
(f_{i}(x),f_{i}(y))\leq \varphi (\delta (x,y))$ \textit{for each }$i\in
\{1,2,...,N\}$\textit{\ and each }$x,y\in K$\textit{.}

\bigskip

\textbf{Definition 2.10 }(iterated function system)\textbf{. }\textit{Given
a complete metric space }$(X,d)$\textit{, an iterated function system is a
pair }$\mathcal{S}=((X,d),(f_{i})_{i\in \{1,2,...,N\}})$\textit{, where }$%
f_{i}:X\rightarrow X$ \textit{is a continuous function\ for each }$i\in
\{1,2,...,N\}$, $N\in \mathbb{N}^{\ast }$\textit{.}

\bigskip

\textbf{3.} \textbf{THE\ RESULTS}

\bigskip

\textbf{Some considerations on} \textbf{iterated function systems consisting
of }$\varphi $\textbf{-contractions}

\bigskip

We start with a result that emphasizes some properties of iterated function
systems consisting of $\varphi $-contractions.

\bigskip

\textbf{Proposition 3.1.} \textit{Let us consider an iterated function system%
} $\mathcal{S}=((X,d),(f_{i})_{i\in I})$ \textit{consisting of} $\varphi $%
\textit{-contractions, where }$\varphi $\textit{\ is a comparison function
and the metric space }$(X,d)$\ \textit{is complete and bounded. Then:}

\textit{a) For every }$\alpha \in \Lambda (I)$\textit{, the set }$\underset{%
n\in \mathbb{N}^{\ast }}{\cap }X_{[\alpha ]_{n}}$\textit{\ has a unique
element which is denoted by }$a_{\alpha }$\textit{.}

\textit{b) If }$a_{\alpha }\neq a_{\beta }$\textit{, where }$\alpha ,\beta
\in \Lambda (I)$\textit{, then there exists }$n_{0}\in \mathbb{N}^{\ast }$%
\textit{\ such that }$X_{[\alpha ]_{n_{0}}}\cap X_{[\beta
]_{n_{0}}}=\emptyset $\textit{.}

\textit{Proof}.

a) Let us consider $\alpha =\alpha _{1}\alpha _{2}...\alpha _{m}...\in
\Lambda (I)$ and $n\in \mathbb{N}^{\ast }$. As for every $x,y\in X_{[\alpha
]_{n}}$ there exist $u,v\in X$ such that $x=f_{\alpha _{1}\alpha
_{2}...\alpha _{n}}(u)$ and $y=f_{\alpha _{1}\alpha _{2}...\alpha _{n}}(v)$,
we have $d(x,y)=d(f_{\alpha _{1}\alpha _{2}...\alpha _{n}}(u),f_{\alpha
_{1}\alpha _{2}...\alpha _{n}}(v))\overset{f_{i}\text{ are }\varphi \text{%
-contractions}}{\leq }\varphi ^{\lbrack n]}(d(u,v))\overset{\varphi \text{
is increasing}}{\leq }\varphi ^{\lbrack n]}(d(X))$, so $d(X_{[\alpha
]_{n}})\leq \varphi ^{\lbrack n]}(d(X))$, hence $d(\overline{X_{[\alpha
]_{n}}})\leq \varphi ^{\lbrack n]}(d(X))$ for every $n\in \mathbb{N}^{\ast }$%
. As $(X,d)$ is complete, making use of Remark 2.2, i) and the fact that $%
\overline{X_{[\alpha ]_{n+1}}}\subseteq \overline{X_{[\alpha ]_{n}}}$ for
every $n\in \mathbb{N}^{\ast }$, we conclude that the set $\underset{n\in 
\mathbb{N}}{\cap }\overline{X_{[\alpha ]_{n}}}$ has one element denoted by $%
a_{\alpha }$, i.e. 
\begin{equation}
\underset{n\in \mathbb{N}}{\cap }\overline{X_{[\alpha ]_{n}}}=\{a_{\alpha }\}%
\text{.}  \tag{1}
\end{equation}%
Let us note that $f_{i}(a_{\alpha })\in f_{i}(\underset{n\in \mathbb{N}%
^{\ast }}{\cap }\overline{X_{[\alpha ]_{n}}})\subseteq \underset{n\in 
\mathbb{N}^{\ast }}{\cap }f_{i}(\overline{X_{[\alpha ]_{n}}})\overset{\text{%
Remark 2.4}}{\subseteq }\underset{n\in \mathbb{N}^{\ast }}{\cap }\overline{%
f_{i}(X_{[\alpha ]_{n}})}=\underset{n\in \mathbb{N}^{\ast }}{\cap }\overline{%
X_{[i\alpha ]_{n}}}\overset{(1)}{=}\{a_{i\alpha }\}$, so 
\begin{equation}
f_{i}(a_{\alpha })=a_{i\alpha }\text{,}  \tag{2}
\end{equation}%
for every $i\in I$ and every $\alpha \in \Lambda (I)$. For $\alpha =\alpha
_{1}\alpha _{2}...\alpha _{n}...\in \Lambda (I)$ and $n\in \mathbb{N}^{\ast
} $, with the notation $\beta _{n}=\alpha _{n+1}\alpha _{n+2}...\alpha
_{m}...\in \Lambda (I)$, we have $a_{\alpha }=a_{[\alpha ]_{n}\beta _{n}}%
\overset{(2)}{=}f_{[\alpha ]_{n}}(a_{\beta _{n}})\in X_{[\alpha ]_{n}}$.
Hence $\{a_{\alpha }\}\subseteq \underset{n\in \mathbb{N}^{\ast }}{\cap }%
X_{[\alpha ]_{n}}\subseteq \underset{n\in \mathbb{N}^{\ast }}{\cap }%
\overline{X_{[\alpha ]_{n}}}=\{a_{\alpha }\}$, so $\{a_{\alpha }\}=\underset{%
n\in \mathbb{N}^{\ast }}{\cap }X_{[\alpha ]_{n}}$.

b) Let us consider $\alpha ,\beta \in \Lambda (I)$ such that\textit{\ }$%
a_{\alpha }\neq a_{\beta }$. Then Remark 2.2, i) assures the existence of a $%
n_{0}\in \mathbb{N}^{\ast }$ such that $\varphi ^{\lbrack n_{0}]}(d(X))<%
\frac{d(a_{\alpha },a_{\beta })}{3}$. Consequently, since (as we have seen
above) $d(X_{[\alpha ]_{n_{0}}})\leq \varphi ^{\lbrack n_{0}]}(d(X))$ and $%
d(X_{[\beta ]_{n_{0}}})\leq \varphi ^{\lbrack n_{0}]}(d(X))$, we get $%
d(X_{[\alpha ]_{n_{0}}})<\frac{d(a_{\alpha },a_{\beta })}{3}$ and $%
d(X_{[\beta ]_{n_{0}}})<\frac{d(a_{\alpha },a_{\beta })}{3}$. If, by
reductio ad absurdum, $X_{[\alpha ]_{n_{0}}}\cap X_{[\beta ]_{n_{0}}}\neq
\emptyset $\textit{, }then choosing $x\in X_{[\alpha ]_{n_{0}}}\cap
X_{[\beta ]_{n_{0}}}$, we get the following contradiction: $d(a_{\alpha
},a_{\beta })\leq d(a_{\alpha },x)+d(x,a_{\beta })\leq d(X_{[\alpha
]_{n_{0}}})+d(X_{[\beta ]_{n_{0}}})<\frac{2d(a_{\alpha },a_{\beta })}{3}$. $%
\square $

\bigskip

\textbf{Remark 3.2.}

i) With the notation\textit{\ }$A=\{a_{\alpha }\mid \alpha \in \Lambda (I)\}$%
, the function $\pi :\Lambda (I)\rightarrow A$, given by $\pi (\alpha
)=a_{\alpha }$ for every $\alpha \in \Lambda (I)$, is continuous.

Indeed, given a fixed $\alpha \in \Lambda (I)$, as $\underset{n\rightarrow
\infty }{\lim }d(X_{[\alpha ]_{n}})=0$, for every $\varepsilon >0$ there
exists $m\in \mathbb{N}^{\ast }$ such that $X_{[\alpha ]_{m}}\subseteq
B(a_{\alpha },\varepsilon )$, so $B(\alpha ,\frac{1}{3^{m}})\subseteq
\{\omega \in \Lambda (I)\mid \lbrack \omega ]_{m}=[\alpha ]_{m}\}\subseteq
\pi ^{-1}(X_{[\alpha ]_{m}})\subseteq \pi ^{-1}(B(a_{\alpha },\varepsilon ))$%
, i.e. $\pi (B(\alpha ,\frac{1}{3^{m}}))\subseteq B(\pi (\alpha
),\varepsilon )$.

ii) Considering the function $F_{\mathcal{S}}:\mathcal{K}(X)\rightarrow 
\mathcal{K}(X)$ given by $F_{\mathcal{S}}(C)=\underset{i\in I}{\cup }%
f_{i}(C) $ for every $C\in \mathcal{K}(X)$, using $(2)$ from the proof of
Proposition 3.1, we infer that $F_{\mathcal{S}}(A)=A$, i.e., taking into
account the uniqueness of the fixed point of $F_{\mathcal{S}}$ (see Theorem
2.5 from [7]), $A$\ is the attractor of the iterated function system $%
\mathcal{S}$. Moreover, the same result guarantees that $\underset{%
n\rightarrow \infty }{\lim }h(F_{\mathcal{S}}^{[n]}(B),A)=0$ for every $B\in 
\mathcal{K}(X)$, where $h$ designates the Hausdorff-Pompeiu metric.

\bigskip

\textbf{The notion of family of functions having attractor}

\bigskip

As $X_{[\alpha ]_{0}}=X$, the above considerations suggest the following:

\bigskip

\textbf{Definition 3.3.} \textit{We say that\ a family of functions }$%
\mathcal{F}=(f_{i})_{i\in I}$\textit{, where }$f_{i}:X\rightarrow X$\textit{%
\ and }$I$\textit{\ is finite, has attractor if the following two properties
are valid:}

\textit{\qquad a) For every }$\alpha \in \Lambda (I)$\textit{, the set }$%
\underset{n\in \mathbb{N}}{\cap }X_{[\alpha ]_{n}}$\textit{\ has a unique
element which is denoted by }$a_{\alpha }$\textit{.}

\textit{\qquad b) If }$a_{\alpha }\neq a_{\beta }$\textit{, where }$\alpha
,\beta \in \Lambda (I)$\textit{, then there exists }$n_{0}\in \mathbb{N}$%
\textit{\ such that }$X_{[\alpha ]_{n_{0}}}\cap X_{[\beta
]_{n_{0}}}=\emptyset $\textit{.}

\textit{The set }$A\overset{def}{=}\{a_{\alpha }\mid \alpha \in \Lambda
(I)\} $ \textit{is called the attractor of} $\mathcal{F}$\textit{.}

\bigskip

\textbf{A metric on the attractor which makes }$\varphi $\textbf{%
-contractions all the functions of a family having attractor}

\bigskip

\textbf{Theorem 3.4.} \textit{If }$\mathcal{F}=(f_{i})_{i\in I}$ \textit{is
a family of functions having attractor }$A$\textit{, then there exist a
metric }$d$\textit{\ on }$A$\textit{\ and a comparison function }$\varphi $%
\textit{\ such that }$d(f_{i}(x),f_{i}(y))\leq \varphi (d(x,y))$\textit{\
for every }$i\in I$ \textit{and every }$x,y\in A$\textit{.}

\textit{Proof}.\textbf{\ }Considering the function $\pi :\Lambda
(I)\rightarrow A$, given by $\pi (\alpha )=a_{\alpha }$ for every $\alpha
\in \Lambda (I)$, the binary relation on $\Lambda (I)$, given by $\alpha
\sim \beta $ if and only if $\pi (\alpha )=\pi (\beta )$, turns out to be an
equivalence relation. We transport the quotient topology on $\Lambda
(I)\diagup \sim $ on the topology $\tau _{A}$ on $A$ via the bijection $%
g:\Lambda (I)\diagup \sim \rightarrow A$ given by $g([\alpha ])=\pi (\alpha
) $ for every $[\alpha ]\in \Lambda (I)\diagup \sim $.

Note that:

i) $g$ is a homeomorphism;

ii) the function $p:\Lambda (I)\rightarrow \Lambda (I)\diagup \sim $, given
by $p(\alpha )=[\alpha ]$ for every $\alpha \in \Lambda (I)$, is continuous;

iii) $\pi =g\circ p$ is continuous.

\medskip

\textit{Claim 1. }$f_{i}\circ \pi =\pi \circ \tau _{i}$ \textit{for every }$%
i\in I$\textit{.}

\textit{Justification of claim 1}. We have $(f_{i}\circ \pi )(\alpha
)=f_{i}(a_{\alpha })\in f_{i}(\underset{n\in \mathbb{N}}{\cap }X_{[\alpha
]_{n}})\subseteq \underset{n\in \mathbb{N}}{\cap }f_{i}(X_{[\alpha ]_{n}})=%
\underset{n\in \mathbb{N}}{\cap }X_{[i\alpha ]_{n}}=\{a_{i\alpha }\}=\{(\pi
\circ \tau _{i})(\alpha )\}$ for every $i\in I$ and every $\alpha \in
\Lambda (I)$.

\bigskip

Note that Claim 1 implies that $A=\underset{i\in I}{\cup }f_{i}(A)$.

\medskip

\textit{Claim 2. }$f_{i}:(A,\tau _{A})\rightarrow (A,\tau _{A})$ \textit{is
continuous for every }$i\in I$\textit{.}

\textit{Justification of claim 2}. Taking into account i), it suffices to
prove that $f_{i}\circ g:\Lambda (I)\diagup \sim \rightarrow A$, given by 
\begin{equation}
(f_{i}\circ g)([\alpha ])\overset{\text{Claim 1}}{=}(\pi \circ \tau
_{i})(\alpha )\text{,}  \tag{1}
\end{equation}%
is continuous. Since $\alpha \sim \beta \Leftrightarrow \pi (\alpha )=\pi
(\beta )\Rightarrow (f_{i}\circ \pi )(\alpha )=(f_{i}\circ \pi )(\beta )%
\overset{\text{Claim 1}}{\Leftrightarrow }(\pi \circ \tau _{i})(\alpha
)=(\pi \circ \tau _{i})(\beta )\overset{(1)}{\Leftrightarrow }(f_{i}\circ
g)([\alpha ])=(f_{i}\circ g)([\beta ])$ and the function $h=\pi \circ \tau
_{i}:\Lambda (I)\rightarrow A$, described by $h(\alpha )=(f_{i}\circ
g)([\alpha ])$ for every $\alpha \in \Lambda (I)$, is continuous (as a
composition of continuous functions; see Remark 2.6 and iii)), relying on
Theorem 4.3, page 126, from [6], we get the conclusion.

\medskip

\textit{Claim 3. }$(A,\tau _{A})$\textit{\ is compact.}

\textit{Justification of claim 3}. From ii) and Remark 2.7, ii), we conclude
that $\Lambda (I)\diagup \sim $ is compact. Using i), we get the conclusion.

\medskip

\textit{Claim 4. The set }$R=\{(\alpha ,\beta )\in \Lambda (I)\times \Lambda
(I)\mid \alpha \sim \beta \}$\textit{\ is closed.}

\textit{Justification of claim 4}. Let us consider $(\alpha ,\beta )\in 
\overline{R}$. Then there exists $((\alpha _{n},\beta _{n}))_{n\in \mathbb{N}%
}\subseteq R$ such that $\underset{n\rightarrow \infty }{\lim }(\alpha
_{n},\beta _{n})=(\alpha ,\beta )$ and consequently $\underset{n\rightarrow
\infty }{\lim }\alpha _{n}=\alpha $ and $\underset{n\rightarrow \infty }{%
\lim }\beta _{n}=\beta $. If $a_{\alpha }\neq a_{\beta }$, then, according
to the property b) from the definition of a family of functions having
attractor, there exists\textit{\ }$n_{0}\in \mathbb{N}$\textit{\ }such that%
\textit{\ }$X_{[\alpha ]_{n_{0}}}\cap X_{[\beta ]_{n_{0}}}=\emptyset $. As $%
\underset{n\rightarrow \infty }{\lim }\alpha _{n}=\alpha $ and $\underset{%
n\rightarrow \infty }{\lim }\beta _{n}=\beta $, there exists $n_{1}\in 
\mathbb{N}$\textit{\ }such that $[\alpha _{n}]_{n_{0}}=[\alpha ]_{n_{0}}$
and $[\beta _{n}]_{n_{0}}=[\beta ]_{n_{0}}$ for every $n\in \mathbb{N}$, $%
n\geq n_{1}$ (see Remark 2.7, i)). But $\alpha _{n}\sim \beta _{n}$ (because 
$(\alpha _{n},\beta _{n})\in R$), i.e. $a_{\alpha _{n}}=a_{\beta _{n}}$, and
therefore we get the following contradiction: $a_{\alpha _{n}}=a_{\beta
_{n}}\in X_{[\alpha _{n}]_{n_{0}}}\cap X_{[\beta _{n}]_{n_{0}}}=X_{[\alpha
]_{n_{0}}}\cap X_{[\beta ]_{n_{0}}}=\emptyset $. Hence $a_{\alpha }=a_{\beta
}$, i.e. $\alpha \sim \beta $, so $(\alpha ,\beta )\in R$. Therefore $R$ is
closed.

\medskip

\textit{Claim 5. }$(A,\tau _{A})$\textit{\ is Hausdorff.}

\textit{Justification of claim 5}. From the compactness of $\Lambda (I)$
(see Remark 2.7, ii)) and Claim 4, we infer that $\Lambda (I)\diagup \sim $
is Hausdorff. Using i) we get the conclusion.

\medskip

Claims 1, 2, 3 and 5 assure us that $(A,(f_{i})_{i\in I})$ is a topological
self-similar system and, based on Theorem 2.9, there exist a metric $d$ on $%
A $ compatible with $\tau _{A}$ and a comparison function\textit{\ }$\varphi
:[0,\infty )\rightarrow \lbrack 0,\infty )$ such that $d(f_{i}(x),f_{i}(y))%
\leq \varphi (d(x),d(y))$ for every $i\in I$ and every $x,y\in A$. $\square $

\bigskip

Let us consider the function $n:X\rightarrow \mathbb{N}\cup \{\infty \}$
given by $n(x)=\sup \{m\in \mathbb{N}\mid x\in F^{[m]}(X)\}$ for every $x\in
X$, where $F:\mathcal{P}(X)\rightarrow \mathcal{P}(X)$ is described by $F(C)=%
\underset{i\in I}{\cup }f_{i}(C)$ for every $C\in \mathcal{P}(X)\overset{def}%
{=}\{Y\mid Y\subseteq X\}$.

\bigskip

The following result provides an alternative characterization of the
attractor $A$ via the function $n$.

\bigskip

\textbf{Proposition 3.5.} \textit{In the framework of the above theorem, we
have }$A=\{x\in X\mid n(x)=\infty \}$.

\textit{Proof}.

"$\subseteq $" If $x\in A$, then there exists $\alpha \in \Lambda (I)$ such
that $x=a_{\alpha }$, hence $x\in X_{[\alpha ]_{m}}\subseteq F^{[m]}(X)$ for
every $m\in \mathbb{N}$. So $n(x)=\sup \{m\in \mathbb{N}\mid x\in
F^{[m]}(X)\}=\sup \mathbb{N}=\infty $.

"$\supseteq $" Since $n(x)=\sup \{m\in \mathbb{N}\mid x\in
F^{[m]}(X)\}=\infty $, for every $m\in \mathbb{N}$ there exists $\alpha
_{m}\in \Lambda _{m}(I)$ such that $x\in X_{\alpha _{m}}$. There exists $%
i_{1}\in I$ such that $\{\omega \in \Lambda ^{\ast }(I)\mid x\in
X_{i_{1}\omega }\}$ is infinite. Indeed, if this is not the case, then the
set $M_{i}\overset{def}{=}\{\omega \in \Lambda ^{\ast }(I)\mid x\in
X_{i\omega }\}$ is finite for every $i\in I$. If $m_{i}\overset{def}{=}\max
\{\left\vert i\omega \right\vert \mid \omega \in M_{i}\}$, then we get the
contradiction that there exists no $\alpha \in \Lambda _{m+1}(I)$ such that $%
x\in X_{\alpha }$, where $m=1+\max \{m_{i}\mid i\in I\}$. Repeating this
procedure we get $\alpha =\alpha _{1}\alpha _{2}...\alpha _{n}...\in \Lambda
(I)$ such that $x\in \underset{n\in \mathbb{N}}{\cap }X_{[\alpha ]_{n}}$,
i.e. $x=a_{\alpha }\in A$. $\square $

\bigskip

\textbf{The family of sets }$\{\overset{\sim }{X_{\alpha }}\mid \alpha \in
\Lambda ^{\ast }(I)\}$\textbf{\ associated to a family of functions} \textbf{%
having attractor}

\bigskip

Given a family of functions\textit{\ }$\mathcal{F}=(f_{i})_{i\in I}$ having
attractor $A$, in the sequel, for $\alpha \in \Lambda ^{\ast }(I)$ we shall
use the following notations:%
\begin{equation*}
Y_{\alpha }\overset{not}{=}\{a_{\beta }\mid X_{\alpha }\cap X_{[\beta
]_{n}}\neq \emptyset \text{ for every }n\in \mathbb{N}\}\text{ and }\overset{%
\sim }{X_{\alpha }}\overset{not}{=}X_{\alpha }\cup Y_{\alpha }\text{.}
\end{equation*}

\bigskip

\textbf{Proposition 3.6 }(The properties of the sets $X_{\alpha }$ and $%
Y_{\alpha }$)\textbf{.} \textit{In the above framework, we have:}

\textit{a)} $A_{\alpha }\subseteq Y_{\alpha }\subseteq A$ \textit{for every} 
$\alpha \in \Lambda ^{\ast }(I)$\textit{;}

\textit{b)} $X_{\alpha }\subseteq \overset{\sim }{X_{\alpha }}\subseteq
X_{\alpha }\cup A$ \textit{for every} $\alpha \in \Lambda ^{\ast }(I)$%
\textit{;}

\textit{c)} $Y_{[\alpha ]_{n+1}}\subseteq Y_{[\alpha ]_{n}}$ \textit{for
every }$\alpha \in \Lambda (I)$ \textit{and every }$n\in \mathbb{N}$\textit{;%
}

\textit{d) }$\underset{n\in \mathbb{N}}{\cap }(X_{[\alpha ]_{n}}\cup
Y_{[\alpha ]_{n}})=(\underset{n\in \mathbb{N}}{\cap }X_{[\alpha ]_{n}})\cup (%
\underset{n\in \mathbb{N}}{\cap }Y_{[\alpha ]_{n}})$ \textit{for every }$%
\alpha \in \Lambda (I)$\textit{;}

\textit{e)} $\underset{n\in \mathbb{N}}{\cap }Y_{[\alpha ]_{n}}=\{a_{\alpha
}\}$ \textit{for every }$\alpha \in \Lambda (I)$\textit{;}

\textit{f)} $A\cap X_{\alpha }\subseteq Y_{\alpha }$ \textit{for every} $%
\alpha \in \Lambda ^{\ast }(I)$\textit{;}

\textit{g)} $f_{i}(Y_{\alpha })\subseteq Y_{i\alpha }$ \textit{for every} $%
\alpha \in \Lambda ^{\ast }(I)$\textit{\ and every }$i\in I$.

\textit{Proof.}

a) If $z\in A_{\alpha }$, then there exists $\gamma \in \Lambda (I)$ such
that $z=f_{\alpha }(a_{\gamma })$, so $z\in X_{\alpha }$. Moreover, $%
z=f_{\alpha }(a_{\gamma })\in f_{\alpha }(\underset{n\in \mathbb{N}}{\cap }%
X_{[\gamma ]_{n}})\subseteq \underset{n\in \mathbb{N}}{\cap }f_{\alpha
}(X_{[\gamma ]_{n}})\subseteq \underset{n\in \mathbb{N}}{\cap }X_{[\alpha
\gamma ]_{n}}=\{a_{\alpha \gamma }\}$, hence $z=a_{\alpha \gamma }\in
X_{\alpha }\cap X_{[\alpha \gamma ]_{n}}$ for every $n\in \mathbb{N}$, i.e. $%
z\in Y_{\alpha }$.

b) It results immediately from a).

c) If $z\in Y_{[\alpha ]_{n+1}}$, then there exists $\beta \in \Lambda (I)$
such that $z=a_{\beta }$ and \linebreak $X_{[\alpha ]_{n+1}}\cap X_{[\beta
]_{k}}\neq \emptyset $ for every $k\in \mathbb{N}$. As $X_{[\alpha
]_{n+1}}\cap X_{[\beta ]_{k}}\subseteq X_{[\alpha ]_{n}}\cap X_{[\beta
]_{k}} $, we deduce that $X_{[\alpha ]_{n}}\cap X_{[\beta ]_{k}}\neq
\emptyset $ for every $k\in \mathbb{N}$, i.e. $z=a_{\beta }\in Y_{[\alpha
]_{n}}$.

d)

"$\supseteq $" It is clear.

"$\subseteq $" Let us suppose that there exists $x\in X_{[\alpha ]_{n}}\cup
Y_{[\alpha ]_{n}}$ for every $n\in \mathbb{N}$ such that $x\notin (\underset{%
n\in \mathbb{N}}{\cap }X_{[\alpha ]_{n}})\cup (\underset{n\in \mathbb{N}}{%
\cap }Y_{[\alpha ]_{n}})$, i.e. there exist $n_{1},n_{2}\in \mathbb{N}$ such
that $x\notin X_{[\alpha ]_{n_{1}}}$ and $x\notin Y_{[\alpha ]_{n_{2}}}$.
Then, in view of c), we have $x\notin X_{[\alpha ]_{m}}$ and $x\notin
Y_{[\alpha ]_{m}}$ which leads to the contradiction $x\notin X_{[\alpha
]_{m}}\cup Y_{[\alpha ]_{m}}$, where $m=\max \{n_{1},n_{2}\}$.

e)

"$\supseteq $" We have $a_{\alpha }\overset{\text{Claim 1 from the proof of
Theorem 3.4}}{\in }\underset{n\in \mathbb{N}}{\cap }A_{[\alpha ]_{n}}\overset%
{\text{a)}}{\subseteq }\underset{n\in \mathbb{N}}{\cap }Y_{[\alpha ]_{n}}$.

"$\subseteq $" If $c\in \underset{n\in \mathbb{N}}{\cap }Y_{[\alpha ]_{n}}$,
then there exists $(\beta _{n})_{n\in \mathbb{N}}\subseteq \pi
^{-1}(\{c\})\subseteq \Lambda (I)$ such that 
\begin{equation}
X_{[\alpha ]_{n}}\cap X_{[\beta _{n}]_{k}}\neq \emptyset \text{,}  \tag{1}
\end{equation}%
for every $n,k\in \mathbb{N}$. The compactness of $\Lambda (I)$ (see Remark
2.7, ii)) assures the existence of a subsequence $(\beta _{n_{l}})_{l\in 
\mathbb{N}}$ of $(\beta _{n})_{n\in \mathbb{N}}$ and of an element $\beta
\in \Lambda (I)$ such that $\underset{l\rightarrow \infty }{\lim }\beta
_{n_{l}}=\beta $. As $\pi (\beta _{n_{l}})=c$, i.e. $a_{\beta _{n_{l}}}=c$,
and $\pi $ is continuous (see Remark 3.2, i)), we infer that $\pi (\beta )=c$%
, i.e. $a_{\beta }=c$. By replacing $\beta _{j}$ with $\beta _{n_{l}}$ for
all $j\in \{n_{l-1}+1,...,n_{l}-1\}$, we can suppose that $\underset{%
n\rightarrow \infty }{\lim }\beta _{n}=\beta $. Hence for every $l\in 
\mathbb{N}$ there exists $n_{l}\in \mathbb{N}$, $n_{l}>l$ such that 
\begin{equation}
\lbrack \beta _{n}]_{l}=[\beta ]_{l}\text{,}  \tag{2}
\end{equation}%
for all $n\in \mathbb{N}$, $n\geq n_{l}$. Hence $X_{[\alpha ]_{n_{l}}}\cap
X_{[\beta _{n_{l}}]_{l}}\overset{(1)}{\neq }\emptyset $, i.e., in view of $%
(2)$, $X_{[\alpha ]_{n_{l}}}\cap X_{[\beta ]_{l}}\neq \emptyset $ and since $%
X_{[\alpha ]_{n_{l}}}\subseteq X_{[\alpha ]_{l}}$, we infer that $X_{[\alpha
]_{l}}\cap X_{[\beta ]_{l}}\neq \emptyset $ for every $l\in \mathbb{N}$.
Therefore, taking into account the property b) of a family of functions
having attractor, we conclude that $a_{\alpha }=a_{\beta }=c$.

f) If $z\in A\cap X_{\alpha }$, then there exists $\gamma \in \Lambda (I)$
such that $z=a_{\gamma }\in X_{[\gamma ]_{n}}$, so $z\in X_{\alpha }\cap
X_{[\gamma ]_{n}}$ and therefore $X_{\alpha }\cap X_{[\gamma ]_{n}}\neq
\emptyset $ for every $n\in \mathbb{N}$. Consequently $z=a_{\gamma }\in
Y_{\alpha }$.

g) If $z\in f_{i}(Y_{\alpha })$, then there exists $\beta \in \Lambda (I)$
such that $z=f_{i}(a_{\beta })$ and $X_{\alpha }\cap X_{[\beta ]_{n}}\neq
\emptyset $ for every $n\in \mathbb{N}$. Since $\emptyset \neq
f_{i}(X_{\alpha }\cap X_{[\beta ]_{n}})\subseteq f_{i}(X_{\alpha })\cap
f_{i}(X_{[\beta ]_{n}})=X_{i\alpha }\cap X_{[i\beta ]_{n+1}}$ for every $%
n\in \mathbb{N}$, we conclude that $a_{i\beta }\overset{\text{Claim 1 from
the proof of Theorem 3.4}}{=}f_{i}(a_{\beta })=z\in Y_{i\alpha }$. $\square $

\bigskip

\textbf{Proposition 3.7} (The properties of the sets $\overset{\sim }{%
X_{\alpha }}$)\textbf{.} \textit{In the above framework, we have:}

\textit{a)} $\overset{\sim }{X_{[\alpha ]_{n+1}}}\overset{\sim }{\subseteq
X_{[\alpha ]_{n}}}$ \textit{for every }$\alpha \in \Lambda (I)$ \textit{and
every }$n\in \mathbb{N}$\textit{;}

\textit{b) }$\underset{n\in \mathbb{N}}{\cap }\overset{\sim }{X_{\alpha }}%
=\{a_{\alpha }\}$ \textit{for every }$\alpha \in \Lambda (I)$.

\textit{c) }$a_{\beta }\in \overset{\sim }{X_{\alpha }}$\textit{, provided
that }$X_{\alpha }\cap \overset{\sim }{X_{[\beta ]_{n}}}\neq \emptyset $%
\textit{\ for every} $n\in \mathbb{N}$\textit{, where }$\alpha \in \Lambda
^{\ast }(I)$\textit{\ and }$\beta \in \Lambda (I)$.

\textit{d) for every }$a_{\alpha },a_{\beta }\in A$\textit{\ such that }$%
a_{\alpha }\neq a_{\beta }$\textit{, there exists} $n_{0}\in \mathbb{N}$ 
\textit{having the property that} $\overset{\sim }{X_{[\alpha ]_{n_{0}}}}%
\cap \overset{\sim }{X_{[\beta ]_{n_{0}}}}=\emptyset $.

\textit{e) }$f_{i}(\overset{\sim }{X_{\alpha }})\subseteq \overset{\sim }{%
X_{i\alpha }}$ \textit{for every }$i\in I$\textit{\ and every }$\alpha \in
\Lambda ^{\ast }(I)$\textit{.}

\textit{Proof.}

a) As $X_{[\alpha ]_{n+1}}\subseteq X_{[\alpha ]_{n}}$ and $Y_{[\alpha
]_{n+1}}\overset{\text{Proposition 3.6, c)}}{\subseteq }Y_{[\alpha ]_{n}}$,
we infer that $X_{[\alpha ]_{n+1}}\cup Y_{[\alpha ]_{n+1}}\subseteq
X_{[\alpha ]_{n}}\cup Y_{[\alpha ]_{n}}$, i.e. $\overset{\sim }{X_{[\alpha
]_{n+1}}}\subseteq \overset{\sim }{X_{[\alpha ]_{n}}}$.

b) We have $\underset{n\in \mathbb{N}}{\cap }\overset{\sim }{X_{[\alpha
]_{n}}}=\underset{n\in \mathbb{N}}{\cap }(X_{[\alpha ]_{n}}\cup Y_{[\alpha
]_{n}})\overset{\text{Proposition 3.6, d) }}{=}(\underset{n\in \mathbb{N}}{%
\cap }X_{[\alpha ]_{n}})\cup (\underset{n\in \mathbb{N}}{\cap }Y_{[\alpha
]_{n}})$\linebreak $\overset{\text{Proposition 3.6, e) }}{=}\{a_{\alpha }\}$.

c) Since $\overset{\sim }{X_{\alpha }}\cap \overset{\sim }{X_{[\beta ]_{l}}}%
\neq \emptyset $\textit{, }i.e. $(X_{\alpha }\cup Y_{\alpha })\cap
(X_{[\beta ]_{l}}\cup Y_{[\beta ]_{l}})\neq \emptyset $\textit{, }we get $%
(X_{\alpha }\cap X_{[\beta ]_{l}})\cup (X_{\alpha }\cap Y_{[\beta
]_{l}})\cup (Y_{\alpha }\cap X_{[\beta ]_{l}})\cup (Y_{\alpha }\cap
Y_{[\beta ]_{l}})\neq \emptyset $ for every $l\in \mathbb{N}$. Thus, at
least one of the sets $\{l\in \mathbb{N}\mid X_{\alpha }\cap X_{[\beta
]_{l}}\neq \emptyset \}$, $\{l\in \mathbb{N}\mid X_{\alpha }\cap Y_{[\beta
]_{l}}\neq \emptyset \}$, $\{l\in \mathbb{N}\mid Y_{\alpha }\cap X_{[\beta
]_{l}}\neq \emptyset \}$ and $\{l\in \mathbb{N}\mid Y_{\alpha }\cap
Y_{[\beta ]_{l}}\neq \emptyset \}$ is infinite.

If $\{l\in \mathbb{N}\mid X_{\alpha }\cap X_{[\beta ]_{l}}\neq \emptyset \}$
is infinite, then, as $X_{[\beta ]_{l+1}}\subseteq X_{[\beta ]_{l}}$ for
every $l\in \mathbb{N}$, we infer that $X_{\alpha }\cap X_{[\beta ]_{l}}\neq
\emptyset $ for every $l\in \mathbb{N}$, so $a_{\beta }\in Y_{\alpha
}\subseteq X_{\alpha }\cup Y_{\alpha }=\overset{\sim }{X_{\alpha }}$.

Since $X_{\alpha }\cap Y_{[\beta ]_{l}}\overset{\text{Proposition 3.6, a)}}{=%
}X_{\alpha }\cap Y_{[\beta ]_{l}}\cap A=(X_{\alpha }\cap A)\cap Y_{[\beta
]_{l}}\overset{\text{Proposition 3.6, f)}}{\subseteq }Y_{\alpha }\cap
Y_{[\beta ]_{l}}$ and $Y_{\alpha }\cap X_{[\beta ]_{l}}\overset{\text{%
Proposition 3.6, a)}}{=}Y_{\alpha }\cap X_{[\beta ]_{l}}\cap A=Y_{\alpha
}\cap (X_{[\beta ]_{l}}\cap A)\overset{\text{Proposition 3.6, f)}}{\subseteq 
}Y_{\alpha }\cap Y_{[\beta ]_{l}}$ for every $l\in \mathbb{N}$, we deduce
that if one of the sets $\{l\in \mathbb{N}\mid X_{\alpha }\cap Y_{[\beta
]_{l}}\neq \emptyset \}$, $\{l\in \mathbb{N}\mid Y_{\alpha }\cap X_{[\beta
]_{l}}\neq \emptyset \}$ and $\{l\in \mathbb{N}\mid Y_{\alpha }\cap
Y_{[\beta ]_{l}}\neq \emptyset \}$ is infinite, then, in view of Proposition
3.6, c), we have%
\begin{equation}
Y_{\alpha }\cap Y_{[\beta ]_{l}}\neq \emptyset \text{,}  \tag{1}
\end{equation}%
for every $l\in \mathbb{N}$.

We are going to prove that $a_{\beta }\in Y_{\alpha }\subseteq X_{\alpha
}\cup Y_{\alpha }=\overset{\sim }{X_{\alpha }}$ and this will closed the
justification of c).

From $(1)$ we deduce that, for every $n\in \mathbb{N}$, there exists $%
a_{n}\in Y_{\alpha }\cap Y_{[\beta ]_{n}}\neq \emptyset $. Consequently we
can find $\gamma _{n},\gamma _{n}^{^{\prime }}\in \Lambda (I)$ such that 
\begin{equation}
a_{n}=a_{\gamma _{n}}=a_{\gamma _{n}^{^{\prime }}}  \tag{2}
\end{equation}%
and 
\begin{equation}
X_{\alpha }\cap X_{[\gamma _{n}]_{l}}\neq \emptyset \text{ and }X_{[\beta
]_{n}}\cap X_{[\gamma _{n}^{^{\prime }}]_{l}}\neq \emptyset \text{,}  \tag{3}
\end{equation}%
for every $l\in \mathbb{N}$. The compactness of $\Lambda (I)$ (see Remark
2.7, i)) assures the existence of the subsequences $(\gamma _{n_{k}})_{k\in 
\mathbb{N}}$ of $(\gamma _{n})_{n\in \mathbb{N}}$ and $(\gamma
_{n_{k}}^{^{\prime }})_{k\in \mathbb{N}}$ of $(\gamma _{n}^{^{\prime
}})_{n\in \mathbb{N}}$ and of the elements $\gamma _{0},\gamma
_{0}^{^{\prime }}\in \Lambda (I)$ such that $\underset{k\rightarrow \infty }{%
\lim }\gamma _{n_{k}}=\gamma _{0}$ and $\underset{k\rightarrow \infty }{\lim 
}\gamma _{n_{k}}^{^{\prime }}=\gamma _{0}^{^{\prime }}$. By replacing $%
\gamma _{n}$ with $\gamma _{n_{k}}$ for all $n\in \{n_{k-1}+1,...,n_{k}-1\}$
and $\gamma _{n}^{^{\prime }}$ with $\gamma _{n_{k}}^{^{\prime }}$ for all $%
n\in \{n_{k-1}+1,...,n_{k}-1\}$, we can suppose that $\underset{n\rightarrow
\infty }{\lim }\gamma _{n}=\gamma _{0}$ and $\underset{n\rightarrow \infty }{%
\lim }\gamma _{n}^{^{\prime }}=\gamma _{0}^{^{\prime }}$. Hence for every $%
l\in \mathbb{N}$ there exists $n_{l}\in \mathbb{N}$, $n_{l}>l$ such that 
\begin{equation}
\lbrack \gamma _{n}]_{l}=[\gamma _{0}]_{l}\text{ and }[\gamma _{n}^{^{\prime
}}]_{l}=[\gamma _{0}^{^{\prime }}]_{l}\text{,}  \tag{4}
\end{equation}%
for every $n\in \mathbb{N}$, $n\geq n_{l}$. Therefore, since $\emptyset 
\overset{\text{(3)}}{\neq }X_{\alpha }\cap X_{[\gamma _{n}]_{l}}\overset{%
\text{(4)}}{=}X_{\alpha }\cap X_{[\gamma _{0}]_{l}}$ for every $l\in \mathbb{%
N}$, we get that 
\begin{equation}
a_{\gamma _{0}}\in Y_{\alpha }\text{.}  \tag{5}
\end{equation}%
Moreover, since $\emptyset \overset{\text{(3)}}{\neq }X_{[\beta
]_{n_{l}}}\cap X_{[\gamma _{n_{l}}^{^{\prime }}]_{l}}\subseteq X_{[\beta
]_{l}}\cap X_{[\gamma _{n_{l}}^{^{\prime }}]_{l}}\overset{\text{(4)}}{=}%
X_{[\beta ]_{l}}\cap X_{[\gamma _{0}^{^{\prime }}]_{l}}$ for every $l\in 
\mathbb{N}$, taking into account the property b) of a family of functions
having attractor, we conclude that%
\begin{equation}
a_{\gamma _{0}^{^{\prime }}}=a_{\beta }\text{.}  \tag{6}
\end{equation}%
Making use of the continuity of $\pi $ we get that $\underset{n\rightarrow
\infty }{\lim }a_{\gamma _{n}}=a_{\gamma _{0}}$ and $\underset{n\rightarrow
\infty }{\lim }a_{\gamma _{n}^{^{\prime }}}=a_{\gamma _{0}^{^{\prime }}}$
and, taking into account $(2)$, we conclude that $a_{\beta }\overset{\text{%
(6)}}{=}a_{\gamma _{0}^{^{\prime }}}=a_{\gamma _{0}}\overset{\text{(5)}}{\in 
}Y_{\alpha }$.

d) If by reductio ad absurdum, we suppose that $\overset{\sim }{X_{[\alpha
]_{n}}}\cap \overset{\sim }{X_{[\beta ]_{n}}}\neq \emptyset $ for every $%
n\in \mathbb{N}$, then we have $\overset{\sim }{\emptyset \neq X_{[\alpha
]_{\max \{k,l\}}}}\cap \overset{\sim }{X_{[\beta ]_{\max \{k,l\}}}}\overset{%
\text{a)}}{\subseteq }\overset{\sim }{X_{[\alpha ]_{k}}}\cap \overset{\sim }{%
X_{[\beta ]_{l}}}$, hence $\overset{\sim }{X_{[\alpha ]_{k}}}\cap \overset{%
\sim }{X_{[\beta ]_{l}}}\neq \emptyset $, for every $k,l\in \mathbb{N}$, so,
based on c), we get that $a_{\beta }\in \overset{\sim }{X_{[\alpha ]_{k}}}$
for every $k\in \mathbb{N}$. Using b) we arrive to the contradiction that $%
a_{\alpha }=a_{\beta }$.

e) We have $f_{i}(\overset{\sim }{X_{\alpha }})=f_{i}(X_{\alpha }\cup
Y_{\alpha })=f_{i}(X_{\alpha })\cup f_{i}(Y_{\alpha })=X_{i\alpha }\cup
f_{i}(Y_{\alpha })$\linebreak $\overset{\text{Proposition 3.6, g)}}{%
\subseteq }X_{i\alpha }\cup Y_{i\alpha }=\overset{\sim }{X_{i\alpha }}$. $%
\square $

\bigskip

\textbf{The semi-metric }$d^{\mu }$\textbf{\ associated to a decreasing
sequence} $\mu $\textbf{\ and to a family of functions} \textbf{having
attractor}

\bigskip

Given a family of functions\textit{\ }$\mathcal{F}=(f_{i})_{i\in I}$ having
attractor and a sequence $\mu =(z_{n})_{n\in \mathbb{N}}$ such that $%
0<z_{n+1}\leq z_{n}$ for every $n\in \mathbb{N}$, we consider the function $%
d^{\mu }:X\times X\rightarrow \lbrack 0,\infty )$ described by $d^{\mu
}(x,y)=\{%
\begin{array}{cc}
0\text{,} & x=y\text{,} \\ 
\inf M_{x,y} & x\neq y%
\end{array}%
$, where $M_{x,y}=\{\overset{n}{\underset{i=0}{\dsum }}z_{\left\vert \alpha
_{i}\right\vert }\mid $there exist $n\in \mathbb{N}$ and $\alpha _{0},\alpha
_{1},...,\alpha _{n}\in \Lambda ^{\ast }(I)$ such that $x\in \overset{\sim }{%
X_{\alpha _{0}}}$, $y\in \overset{\sim }{X_{\alpha _{n}}}$ and $\overset{%
\sim }{X_{\alpha _{i}}}\cap \overset{\sim }{X_{\alpha _{i+1}}}\neq \emptyset 
$ for every $i\in \{0,1,...,n-1\}\}$.

\bigskip

\textbf{Proposition 3.8} (The properties of $d^{\mu }$)\textbf{.} \textit{In
the above framework, we have:}

\textit{a)} $d^{\mu }(x,x)=0$ \textit{for every} $x\in X$\textit{;}

\textit{b)} $d^{\mu }(x,y)=d^{\mu }(y,x)$ \textit{for every} $x,y\in X$%
\textit{;}

\textit{c)} $d^{\mu }(x,y)\leq d^{\mu }(x,z)+d^{\mu }(z,y)$ \textit{for every%
} $x,y,z\in X$\textit{;}

\textit{d)} $d^{\mu }(x,y)>0$ \textit{for every }$x\in X\smallsetminus A$%
\textit{\ and every }$y\in X\smallsetminus \{x\}$\textit{;}

\textit{e)} $d^{\mu }(f_{i}(x),f_{i}(y))\leq d^{\mu }(x,y)$ \textit{for every%
} $x,y\in X$\textit{;}

\textit{f)} $d^{\mu }(x,y)\leq z_{0}$ \textit{for every} $x,y\in X$\textit{;}

\textit{g) If the sequence }$\mu $\textit{\ is constant, then }$d^{\mu }$%
\textit{\ is a metric.}

\textit{Proof}.

a) and b) are obvious, while c) is clear since every chain from $x$ to $z$
and every chain from $z$ to $y$ generate a chain from $x$ to $y$.

d) Considering the function $m:X\rightarrow \mathbb{N}\cup \{\infty \}$,
given by $m(u)=\sup \{\left\vert \alpha \right\vert \mid \alpha \in \Lambda
^{\ast }(I)$ and $u\in \overset{\sim }{X_{\alpha }}\}$ for every $u\in X$,
using a similar argument as in the one used in the proof of Proposition 3.5,
we obtain that $m(u)=\infty $ if and only if $u\in A$. Hence $m(x)\in 
\mathbb{N}$ since $x\in X\smallsetminus A$ and if for $n\in \mathbb{N}$ and $%
\alpha _{0},\alpha _{1},...,\alpha _{n}\in \Lambda ^{\ast }(I)$ we have $%
x\in \overset{\sim }{X_{\alpha _{0}}}$, $y\in \overset{\sim }{X_{\alpha _{n}}%
}$ and $\overset{\sim }{X_{\alpha _{i}}}\cap \overset{\sim }{X_{\alpha
_{i+1}}}\neq \emptyset $ for every $i\in \{0,...,n-1\}$, then $z_{m(x)}\leq
z_{\left\vert \alpha _{0}\right\vert }\leq \overset{n}{\underset{i=0}{\dsum }%
}z_{\left\vert \alpha _{i}\right\vert }$. So $z_{m(x)}$ is a lower bound for 
$M_{x,y}$ and, consequently, $0<z_{m(x)}\leq \inf M_{x,y}=d^{\mu }(x,y)$.

e) The inequality is obvious if $f_{i}(x)=f_{i}(y)$ (in particular, if $x=y$%
). Otherwise, if for $n\in \mathbb{N}$ and $\alpha _{0},\alpha
_{1},...,\alpha _{n}\in \Lambda ^{\ast }(I)$ we have $x\in \overset{\sim }{%
X_{\alpha _{0}}}$, $y\in \overset{\sim }{X_{\alpha _{n}}}$ and $\overset{%
\sim }{X_{\alpha _{j}}}\cap \overset{\sim }{X_{\alpha _{j+1}}}\neq \emptyset 
$ for every $j\in \{0,1,...,n-1\}$, then $f_{i}(x)\in f_{i}(\overset{\sim }{%
X_{\alpha _{0}}})\overset{\text{Proposition 3.7, e)}}{\subseteq }\overset{%
\sim }{X_{i\alpha _{0}}}$, $f_{i}(y)\in f_{i}(\overset{\sim }{X_{\alpha _{n}}%
})\overset{\text{Proposition 3.7, e)}}{\subseteq }\overset{\sim }{X_{i\alpha
_{n}}}$ and $\emptyset \neq f_{i}(\overset{\sim }{X_{\alpha _{j}}}\cap 
\overset{\sim }{X_{\alpha _{j+1}}})\subseteq f_{i}(\overset{\sim }{X_{\alpha
_{j}}})\cap f_{i}(\overset{\sim }{X_{\alpha _{j+1}}})\overset{\text{%
Proposition 3.7, e)}}{\subseteq }\overset{\sim }{X_{i\alpha _{j}}}\cap 
\overset{\sim }{X_{i\alpha _{j+1}}}$, so\linebreak\ $\overset{\sim }{%
X_{i\alpha _{j}}}\cap \overset{\sim }{X_{i\alpha _{j+1}}}\neq \emptyset $
for every $j\in \{0,1,...,n-1\}$, i.e. $\overset{n}{\underset{j=0}{\dsum }}%
z_{\left\vert i\alpha _{j}\right\vert }\in M_{f_{i}(x),f_{i}(y)}$. Hence $%
d^{\mu }(f_{i}(x),f_{i}(y))=\inf M_{f_{i}(x),f_{i}(y)}\leq \overset{n}{%
\underset{j=0}{\dsum }}z_{\left\vert i\alpha _{j}\right\vert }\leq \overset{n%
}{\underset{j=0}{\dsum }}z_{\left\vert \alpha _{j}\right\vert }$, i.e. $%
d^{\mu }(f_{i}(x),f_{i}(y))$ is a lower bound for $M_{x,y}$, so $d^{\mu
}(f_{i}(x),f_{i}(y))\leq \inf M_{x,y}=d^{\mu }(x,y)$.

f) and g) result from the definition of $d^{\mu }$. $\square $

\bigskip

\textbf{Remark 3.9}.

i) $d^{\mu }$ is a semi-metric.

ii) From the proof of d) we get that $\{y\in X\mid d^{\mu }(x,y)<\frac{%
z_{m(x)}}{2}\}=\{x\}$ for every $x\in X\smallsetminus A$. In other words,
the topology generated by $d^{\mu }$ on $X\smallsetminus A$ is the discrete
one.

iii) From e) we conclude that (with respect to $d^{\mu }$) each of the
functions $f_{i}$ has the Lipschitz constant less or equal to $1$.

\bigskip

Given a natural number $N$, a family of functions\textit{\ }$\mathcal{F}%
=(f_{i})_{i\in I}$ having attractor and a sequence $\mu =(z_{n})_{n\in 
\mathbb{N}}$ such that $0<z_{n+1}\leq z_{n}$ for every $n\in \mathbb{N}$, we
consider the function $d_{N}^{\mu }:X\times X\rightarrow \lbrack 0,\infty )$
described by $d_{N}^{\mu }(x,y)=\{%
\begin{array}{cc}
0\text{,} & x=y \\ 
\inf M_{x,y}^{N} & x\neq y%
\end{array}%
$, where $M_{x,y}^{N}=\{\overset{n}{\underset{i=0}{\dsum }}z_{\left\vert
\alpha _{i}\right\vert }\mid $there exist $n\in \mathbb{N}$ and $\alpha
_{0},\alpha _{1},...,\alpha _{n}\in \Lambda _{0}(I)\cup \Lambda _{1}(I)\cup
...\cup \Lambda _{N}(I)$ such that $x\in \overset{\sim }{X_{\alpha _{0}}}$, $%
y\in \overset{\sim }{X_{\alpha _{n}}}$ and $\overset{\sim }{X_{\alpha _{i}}}%
\cap \overset{\sim }{X_{\alpha _{i+1}}}\neq \emptyset $ for every $i\in
\{0,1,...,n-1\}\}$.

We also consider the sequences $\mu _{N}$ and $\mu _{N,p}$ given by $\mu
_{N}=(z_{n}^{N})_{n\in \mathbb{N}}$ and $\mu _{N,p}=(z_{n}^{N,p})_{n\in 
\mathbb{N}}$, where $p\in \mathbb{N}$, $z_{n}^{N}=\{%
\begin{array}{cc}
z_{n}\text{,} & n\in \{0,1,...,N\} \\ 
z_{N}\text{,} & n\in \mathbb{N}\text{, }n\geq N+1%
\end{array}%
$ and $z_{n}^{N,p}=\{%
\begin{array}{cc}
z_{n}\text{,} & n\in \{0,1,...,N\} \\ 
z_{N}\text{,} & n\in \{N+1,...,N+p\} \\ 
\frac{z_{N}}{2}\text{,} & n\in \mathbb{N}\text{, }n\geq N+p+1%
\end{array}%
$.

\bigskip

\textbf{Proposition 3.10} (The properties of $d_{N}^{\mu }$)\textbf{.} 
\textit{In the above framework, we have:}

\textit{a)} $d^{\mu }\leq d_{N+1}^{\mu }\leq d_{N}^{\mu }$ \textit{for every 
}$N\in \mathbb{N}$\textit{;}

\textit{b)} $d^{\mu }=\underset{N\rightarrow \infty }{\lim }d_{N}^{\mu }=%
\underset{N\in \mathbb{N}}{\inf }d_{N}^{\mu }$ \textit{(i.e. }$\underset{%
N\rightarrow \infty }{\lim }d_{N}^{\mu }(x,y)=d^{\mu }(x,y)\ $\textit{for
every }$x,y\in X$\textit{);}

\textit{c)} $d^{\mu _{N,p}}\leq d^{\mu _{N,p+1}}\leq d^{\mu _{N}}$ \textit{%
for every }$N,p\in \mathbb{N}$\textit{;}

\textit{d)} $d_{N}^{\mu }=d^{\mu _{N}}$ \textit{for every }$N\in \mathbb{N}$.

\textit{Proof}.

a) It results from the inclusion $M_{x,y}^{N}\subseteq
M_{x,y}^{N+1}\subseteq M_{x,y}$ which is valid for every $x,y\in X$ and
every $N\in \mathbb{N}$.

b) Given $x,y\in X$, $x\neq y$, for every $\varepsilon >0$ there exist $n\in 
\mathbb{N}$ and $\alpha _{0},\alpha _{1},...,\alpha _{n}\in \Lambda ^{\ast
}(I)$ such that $x\in \overset{\sim }{X_{\alpha _{0}}}$, $y\in \overset{\sim 
}{X_{\alpha _{n}}}$ and $\overset{\sim }{X_{\alpha _{i}}}\cap \overset{\sim }%
{X_{\alpha _{i+1}}}\neq \emptyset $ for every $i\in \{0,1,...,n-1\}$ and $%
\overset{n}{\underset{i=0}{\dsum }}z_{\left\vert \alpha _{i}\right\vert
}<d^{\mu }(x,y)+\varepsilon $. As $\overset{n}{\underset{i=0}{\dsum }}%
z_{\left\vert \alpha _{i}\right\vert }\in M_{x,y}^{N_{\varepsilon }}$, where 
$N_{\varepsilon }=\max \{\left\vert \alpha _{0}\right\vert ,\left\vert
\alpha _{1}\right\vert ,...,\left\vert \alpha _{n}\right\vert \}$, we infer
that $d_{N}^{\mu }\leq d_{N_{\varepsilon }}^{\mu }=\inf
M_{x,y}^{N_{\varepsilon }}\leq \overset{n}{\underset{i=0}{\dsum }}%
z_{\left\vert \alpha _{i}\right\vert }<d^{\mu }(x,y)+\varepsilon $, so $%
\left\vert d_{N}^{\mu }-d^{\mu }\right\vert <\varepsilon $ for every $N\in 
\mathbb{N}$, $N\geq N_{\varepsilon }$. Hence $\underset{N\rightarrow \infty }%
{\lim }d_{N}^{\mu }(x,y)=d^{\mu }(x,y)\ $for every $x,y\in X$, $x\neq y$.
The equality is obvious for $x=y$.

c) Let us note that if the sequences $\mu =(z_{n})_{n\in \mathbb{N}}$ and $%
\nu =(t_{n})_{n\in \mathbb{N}}$ have the property $z_{n}\leq t_{n}$ for
every $n\in \mathbb{N}$ (we denote this situation by $\mu \prec \nu $), then 
$d^{\mu }\leq d^{\nu }$. Now the conclusion results from the fact that $\mu
_{N,p}\prec \mu _{N,p+1}\prec \mu _{N}$.

d) First let us note that%
\begin{equation}
d_{N}^{\mu }=d_{N}^{\mu _{N}}\text{.}  \tag{1}
\end{equation}

Moreover%
\begin{equation}
d_{N}^{\mu _{N}}=d_{M}^{\mu _{N}}\text{,}  \tag{2}
\end{equation}%
for every $M\in \mathbb{N}$, $M>N$.

Indeed, we have $d_{M}^{\mu _{N}}\overset{a)}{\leq }d_{N}^{\mu _{N}}$ for
every $M\in \mathbb{N}$, $M>N$, it remains to prove that $d_{N}^{\mu
_{N}}\leq d_{M}^{\mu _{N}}$ for every $M\in \mathbb{N}$, $M>N$. This follows
from the fact that if for $x,y\in X$, $x\neq y$ there exist $n\in \mathbb{N}$
and $\alpha _{0},\alpha _{1},...,\alpha _{n}\in \Lambda _{0}(I)\cup \Lambda
_{1}(I)\cup ...\cup \Lambda _{M}(I)$ such that $x\in \overset{\sim }{%
X_{\alpha _{0}}}$, $y\in \overset{\sim }{X_{\alpha _{n}}}$ and $\overset{%
\sim }{X_{\alpha _{i}}}\cap \overset{\sim }{X_{\alpha _{i+1}}}\neq \emptyset 
$ for every $i\in \{0,1,...,n-1\}$, then $x\in \overset{\sim }{%
X_{^{N}[\alpha _{0}]}}$, $y\in \overset{\sim }{X_{^{N}[\alpha _{n}]}}$, $%
\emptyset \neq \overset{\sim }{X_{\alpha _{i}}}\cap \overset{\sim }{%
X_{\alpha _{i+1}}}\overset{\sim }{\overset{\text{Proposition 3.7, a)}}{%
\subseteq }\overset{\sim }{X_{^{N}[\alpha _{i}]}}}\cap \overset{\sim }{%
X_{_{^{N}[\alpha _{i+1}]}}}$ (so $\overset{\sim }{X_{^{N}[\alpha _{i}]}}\cap 
\overset{\sim }{X_{_{^{N}[\alpha _{i+1}]}}}\neq \emptyset $) for every $i\in
\{0,1,...,n-1\}$ and $\overset{n}{\underset{i=0}{\dsum }}z_{\left\vert
\alpha _{i}\right\vert }=\overset{n}{\underset{i=0}{\dsum }}z_{\left\vert
^{N}[\alpha _{i}]\right\vert }$, where $^{N}[\alpha ]=\{%
\begin{array}{cc}
\alpha \text{,} & \text{if }\left\vert \alpha \right\vert \leq N \\ 
\lbrack \alpha ]_{N}\text{,} & \text{if }\left\vert \alpha \right\vert >N%
\end{array}%
$.

Based on b), by passing to limit as $M$ goes to $\infty $ in\ $(2)$, and
using $(1)$, we get the conclusion. $\square $

\bigskip

\textbf{Proposition 3.11.} \textit{In the above framework, we have} $%
\underset{p\rightarrow \infty }{\lim }d^{\mu _{N,p}}=d^{\mu _{N}}$ \textit{%
(i.e. }$\underset{p\rightarrow \infty }{\lim }d^{\mu _{N,p}}(x,y)=d^{\mu
_{N}}(x,y)$\textit{\ for every }$x,y\in X$\textit{\ ) for every} $N\in 
\mathbb{N}$\textit{.}

\textit{Proof}. Note that $\underset{p\rightarrow \infty }{\lim }d^{\mu
_{N,p}}(x,y)$ exists and is finite since, according to Proposition 3.10, c)
the sequence $(d^{\mu _{N,p}}(x,y))_{p\in \mathbb{N}}$ is increasing and
bounded for every $x,y\in X$.

If $d^{\mu _{N}}(x,y)=0$, then, based on Proposition 3.10, c), we get
that\linebreak $\underset{p\rightarrow \infty }{\lim }d^{\mu
_{N,p}}(x,y)=d^{\mu _{N}}(x,y)$.

Hence we have to consider only the case when $d^{\mu _{N}}(x,y)\neq 0$.
Taking into account Proposition 3.10, c), we have $\underset{p\rightarrow
\infty }{\lim }d^{\mu _{N,p}}(x,y)\leq d^{\mu _{N}}(x,y)$. Let us suppose,
by reductio ad absurdum, that $l_{0}\overset{not}{=}\underset{p\rightarrow
\infty }{\lim }d^{\mu _{N,p}}(x,y)=\underset{p\in \mathbb{N}}{\sup }d^{\mu
_{N,p}}(x,y)<d^{\mu _{N}}(x,y)$. Then $d^{\mu _{N,p}}(x,y)\leq l_{0}<l%
\overset{not}{=}\frac{l_{0}+d^{\mu _{N}}(x,y)}{2}<d^{\mu _{N}}(x,y)$ for
every $p\in \mathbb{N}$. Hence there exist $n_{p}\in \mathbb{N}$ and $\alpha
_{0}^{p},\alpha _{1}^{p},...,\alpha _{n_{p}}^{p}\in \Lambda ^{\ast }(I)$
such that $x\in \overset{\sim }{X_{\alpha _{0}^{p}}}$, $y\in \overset{\sim }{%
X_{\alpha _{n_{p}}^{p}}}$, $\overset{\sim }{X_{\alpha _{i}^{p}}}\cap \overset%
{\sim }{X_{\alpha _{i+1}^{p}}}\neq \emptyset $ for every $i\in
\{0,1,...,n_{p}-1\}$ and $\overset{n_{p}}{\underset{i=0}{\dsum }}%
z_{\left\vert \alpha _{i}^{p}\right\vert }^{N,p}<l$. Since $z_{\left\vert
\alpha _{i}^{p}\right\vert }^{N,p}\geq \frac{z_{N}}{2}$ for every $i\in
\{0,1,...,n_{p}\}$, we infer that $(n_{p}+1)\frac{z_{N}}{2}<l$ for every $%
p\in \mathbb{N}$, so the sequence $(n_{p})_{p\in \mathbb{N}}\subseteq 
\mathbb{N}$ is bounded and therefore there exists a subsequence $%
(n_{p_{k}})_{k\in \mathbb{N}}$ of $(n_{p})_{p\in \mathbb{N}}$ such that $%
n_{p_{1}}=n_{p_{2}}=...\overset{not}{=}m$.

We say that $i\in \{0,1,...,m\}$ is:

- of type $I$ if $\overline{\underset{k\rightarrow \infty }{\lim }}%
\left\vert \alpha _{i}^{p_{k}}\right\vert <\infty $;

- of type $II$ if $\overline{\underset{k\rightarrow \infty }{\lim }}%
\left\vert \alpha _{i}^{p_{k}}\right\vert =\infty $.

If $i$ is of type $I$, then there exists $C\in \mathbb{R}$ such that $%
\left\vert \alpha _{i}^{p_{k}}\right\vert <C$ for every $k\in \mathbb{N}$,
so, by passing to a subsequence, we can assume that $\alpha
_{i}^{p_{1}}=\alpha _{i}^{p_{2}}=...=\alpha _{i}^{p_{n}}=...\overset{not}{=}%
\alpha _{i}$.

If $i$ is of type $II$, then, by passing to a subsequence, we can assume
that:

i) $\underset{k\rightarrow \infty }{\lim }\left\vert \alpha
_{i}^{p_{k}}\right\vert =\infty $;

ii) $\left\vert \alpha _{i}^{p_{k}}\right\vert <\left\vert \alpha
_{i}^{p_{k+1}}\right\vert $ for every $k\in \mathbb{N}$;

iii) there exists $\alpha _{i}\in \Lambda (I)$ such that $[\alpha
_{i}]_{\left\vert \alpha _{i}^{p_{k}}\right\vert }=\alpha _{i}^{p_{k}}$ for
every $k\in \mathbb{N}$ (since there exists $j_{1}\in I$ which is the first
letter for an infinity of $\alpha _{i}^{p_{k}}$ -otherwise, we contradict
i)- and we choose $j_{1}$ to be the first letter of $\alpha _{i}$; the same
argument provides $j_{2}\in I$ which is the second letter for an infinity of 
$\alpha _{i}^{p_{k}}$ having $j_{1}$ as the first letter and we choose $%
j_{2} $ to be the second letter of $\alpha _{i}$; we continue this
procedure).

For a fixed $j\in \{0,1,...,m-1\}$ the following four cases are possible:

a) $j$ and $j+1$ are of type $I$;

b) $j$ is of type $I$ and $j+1$ is of type $II$;

c) $j$ is of type $II$ and $j+1$ is of type $I$;

d) $j$ and $j+1$ are of type $II$.

In case a) we have 
\begin{equation}
\overset{\sim }{X_{\alpha _{j}}}\cap \overset{\sim }{X_{\alpha _{j+1}}}\neq
\emptyset \text{.}  \tag{1}
\end{equation}

In case b) we have $\overset{\sim }{X_{\alpha _{j}}}\cap \overset{\sim }{%
X_{[\alpha _{j+1}]_{\left\vert \alpha _{j+1}^{p_{k}}\right\vert }}}\neq
\emptyset $ for every $k\in \mathbb{N}$, so, according to Proposition 3.7,
c), we get 
\begin{equation}
a_{\alpha _{j+1}}\in \overset{\sim }{X_{\alpha _{j}}}\text{.}  \tag{2}
\end{equation}

In case c), as above, we get 
\begin{equation}
a_{\alpha _{j}}\in \overset{\sim }{X_{\alpha _{j+1}}}\text{.}  \tag{3}
\end{equation}

In case d), we have $\overset{\sim }{X_{[\alpha _{j}]_{\left\vert \alpha
_{j}^{p_{k}}\right\vert }}}\cap \overset{\sim }{X_{[\alpha
_{j+1}]_{\left\vert \alpha _{j+1}^{p_{k}}\right\vert }}}\neq \emptyset $ for
every $k\in \mathbb{N}$, so, using Proposition 3.7, d), we obtain that 
\begin{equation}
a_{\alpha _{j}}=a_{\alpha _{j+1}}\text{.}  \tag{4}
\end{equation}

First let us note that if all $i\in \{0,1,...,m\}$ would be of type $II$,
then $x\in \overset{\sim }{X_{\alpha _{0}^{p_{k}}}}=\overset{\sim }{%
X_{[\alpha _{0}]_{\left\vert \alpha _{0}^{p_{k}}\right\vert }}}$ for every $%
k\in \mathbb{N}$, so taking into account Proposition 3.7, b), we get that $%
x=a_{\alpha _{0}}$. In the same way we obtain that $y=a_{\alpha _{m}}$ and,
based on $(4)$, we conclude that $x=a_{\alpha _{0}}=a_{\alpha
_{1}}=...=a_{\alpha _{m-1}}=a_{\alpha _{m}}=y$ which contradicts our
assumption that $d^{\mu _{N}}(x,y)\neq 0$. Hence we can assume that at least
one $i\in \{0,1,...,m\}$ is of type $I$.

Now we mention the following four facts:

\textit{Fact 1}. As we have seen before, if $i\in \{0,1,...,m\}$ is of type $%
II$ and $u\in \overset{\sim }{X_{\alpha _{i}^{p_{k}}}}$ for every $k\in 
\mathbb{N}$, then $u=a_{\alpha _{i}}$.

\textit{Fact 2}. If $j$ and $q$ are of type $I$ and $j+1$, $j+2$, ..., $q-1$
are of type $II$, where $0\leq j<q\leq m$, then, $a_{\alpha _{j+1}}\overset{%
\text{(2)}}{\in }\overset{\sim }{X_{\alpha _{j}}}$, $a_{\alpha _{q-1}}%
\overset{\text{(3)}}{\in }\overset{\sim }{X_{\alpha _{q}}}$ and, based on $%
(4)$, we have $a_{\alpha _{j+1}}=a_{\alpha _{j+2}}=...=a_{\alpha _{q-1}}$,
so $\overset{\sim }{X_{\alpha _{j}}}\cap \overset{\sim }{X_{\alpha _{q}}}%
\neq \emptyset $.

\textit{Fact 3}. If $j$, $j+1$, ..., $q-1$ are of type $II$ and $q$ is of
type $I$, where $0\leq j<q\leq m$, then based on $(4)$, we have $a_{\alpha
_{j}}=a_{\alpha _{j+1}}=...=a_{\alpha _{q-1}}$ and $a_{\alpha _{q-1}}\overset%
{\text{(3)}}{\in }\overset{\sim }{X_{\alpha _{q}}}$, so $a_{\alpha _{j}}\in 
\overset{\sim }{X_{\alpha _{q}}}$.

\textit{Fact 4}. If $j$ is of type $I$ and $j+1$, ..., $q$ are of type $II$,
where $0\leq j<q\leq m$, then based on $(4)$, we have $a_{\alpha
_{j+1}}=a_{\alpha _{j+2}}=...=a_{\alpha _{q}}$ and $a_{\alpha _{j+1}}\overset%
{\text{(3)}}{\in }\overset{\sim }{X_{\alpha _{j}}}$, so $a_{\alpha _{q}}\in 
\overset{\sim }{X_{\alpha _{j}}}$.

In view of the above mentioned four facts, we can pick up from the set $%
\{\alpha _{0},\alpha _{1},...,\alpha _{m}\}$ a subset $\{\alpha _{i_{0}}%
\overset{not}{=}\beta _{0},\alpha _{i_{1}}\overset{not}{=}\beta
_{1},...,\alpha _{i_{l}}\overset{not}{=}\beta _{l}\}$, where $%
i_{0},i_{1},...,i_{l}$ are all the type $I$ elements of $\{0,1,...,m\}$,
such that $x\in \overset{\sim }{X_{\beta _{0}}}$, $y\in \overset{\sim }{%
X_{\beta _{l}}}$ and $\overset{\sim }{X_{\beta _{j}}}\cap \overset{\sim }{%
X_{\beta _{j+1}}}\neq \emptyset $ for every $j\in \{0,1,...,l-1\}$. Then we
get the following contradiction: $d^{\mu _{N}}(x,y)\leq \overset{l}{\underset%
{j=0}{\dsum }}z_{\left\vert \beta _{j}\right\vert }^{N}\leq \overset{%
n_{p_{k}}}{\underset{i=0}{\dsum }}z_{\left\vert \alpha
_{i}^{p_{k}}\right\vert }^{N}=\overset{n_{p_{k}}}{\underset{i=0}{\dsum }}%
z_{\left\vert \alpha _{i}^{p_{k}}\right\vert }^{N,n_{p_{k}}}<l<d^{\mu
_{N}}(x,y)$, where $k$ is chosen such that $N+p_{k}>\max \{\left\vert \alpha
_{0}\right\vert ,\left\vert \alpha _{1}\right\vert ,...,\left\vert \alpha
_{m}\right\vert \}$. $\square $

\bigskip

\textbf{Proposition 3.12.} \textit{In the above framework, for every }$%
x,y\in X$, $x\neq y$ \textit{and} $M>0$\textit{, there exists a decreasing
sequence} $\mu =(z_{n})_{n\in \mathbb{N}}$ \textit{such that:}

\textit{i)} $\underset{n\rightarrow \infty }{\lim }z_{n}=0$;

\textit{ii)} $d^{\mu }(x,y)>0$\textit{;}

\textit{iii)} $d^{\mu }\leq M$.

\textit{Proof}. For $M>0$, let us consider the sequence $\mu ^{0}$, where $%
\mu ^{0}=(y_{n})_{n\in \mathbb{N}}$ and $y_{n}=M$ for every $n\in \mathbb{N}$%
. Note that $d^{\mu ^{0}}(x,y)=M$. By mathematical induction we construct a
sequence $(\mu ^{k})_{k\in \mathbb{N}}$ of sequences such that%
\begin{equation}
d^{\mu ^{k}}(x,y)-d^{\mu ^{k+1}}(x,y)<\frac{M}{2^{k+2}}\text{,}  \tag{1}
\end{equation}%
for every $k\in \mathbb{N}$. In fact we construct a strictly increasing
sequence $(p_{k})_{k\in \mathbb{N}}\subseteq \mathbb{N}$ such that $\mu
^{k+1}=\mu _{p_{k}+1,p_{k+1}-p_{k}}^{k}$, where if $p_{k}$ is constructed, $%
p_{k+1}$ is chosen such that $(1)$ is valid based on the fact that $\underset%
{p\rightarrow \infty }{\lim }d^{\mu _{p_{k}+1,p}^{k}}=d^{\mu _{p_{k}+1}^{k}}$
(see Proposition 3.11). Note that 
\begin{equation}
\mu ^{k+1}\prec \mu ^{k}\text{,}  \tag{2}
\end{equation}%
for every $k\in \mathbb{N}$. Since $\left\Vert \mu ^{k+1}-\mu
^{k}\right\Vert \leq \frac{M}{2^{k+1}}$ for every $k\in \mathbb{N}$ (here,
for a sequence $(a_{n})_{n\in \mathbb{N}}$, by $\left\Vert (a_{n})_{n\in 
\mathbb{N}}\right\Vert $ we mean $\underset{n\in \mathbb{N}}{\sup }%
\left\vert a_{n}\right\vert $), we infer that the sequence $(\mu ^{k})_{k\in 
\mathbb{N}}$ is Cauchy, so it is convergent and therefore there exists a
sequence $\mu =(z_{n})_{n\in \mathbb{N}}$ such that $\mu =\underset{%
k\rightarrow \infty }{\lim }\mu ^{k}$.

Let us note that we have $z_{n}=\{%
\begin{array}{cc}
M\text{,} & n\in \{0,1,...,p_{1}\} \\ 
\frac{M}{2}\text{,} & n\in \{p_{1}+1,...,p_{2}\} \\ 
\frac{M}{2^{2}}\text{,} & n\in \{p_{2}+1,...,p_{3}\} \\ 
... & ... \\ 
\frac{M}{2^{q}}\text{,} & n\in \{p_{q}+1,...,p_{q+1}\} \\ 
... & ...%
\end{array}%
$ and, if \linebreak $\mu ^{k}=(z_{n}^{k})_{n\in \mathbb{N}}$, then $%
z_{n}^{k}=\{%
\begin{array}{cc}
M\text{,} & n\in \{0,1,...,p_{1}\} \\ 
\frac{M}{2}\text{,} & n\in \{p_{1}+1,...,p_{2}\} \\ 
... & ... \\ 
\frac{M}{2^{k-1}}\text{,} & n\in \{p_{k-1}+1,...,p_{k}\} \\ 
\frac{M}{2^{k}}\text{,} & n\geq p_{k}+1%
\end{array}%
$ for every $k\in \mathbb{N}$.

Now we prove that the decreasing sequence $\mu $ satisfies the conditions
i), ii) and iii).

i) is obvious having in view the description of the general term of $\mu $.

ii) First of all let us note that 
\begin{equation}
d^{\mu }(x_{1},y_{1})=\underset{k\rightarrow \infty }{\lim }d^{\mu
^{k}}(x_{1},y_{1})=\underset{k\in \mathbb{N}}{\inf }d^{\mu ^{k}}(x_{1},y_{1})%
\text{,}  \tag{3}
\end{equation}%
for every $x_{1},y_{1}\in X$.

Indeed, let us fix $x_{1},y_{1}\in X$. For every $\varepsilon >0$ there
exist $n\in \mathbb{N}$ and $\alpha _{0},\alpha _{1},...,\alpha _{n}\in
\Lambda ^{\ast }(I)$ such that $x_{1}\in \overset{\sim }{X_{\alpha _{0}}}$, $%
y_{1}\in \overset{\sim }{X_{\alpha _{n}}}$, $\overset{\sim }{X_{\alpha _{i}}}%
\cap \overset{\sim }{X_{\alpha _{i+1}}}\neq \emptyset $ for every $i\in
\{0,1,...,n-1\}$ and $\overset{n}{\underset{i=0}{\dsum }}z_{\left\vert
\alpha _{i}\right\vert }<d^{\mu }(x_{1},y_{1})+\varepsilon $. There exists $%
k_{\varepsilon }\in \mathbb{N}$ such that $z_{\left\vert \alpha
_{i}\right\vert }=z_{\left\vert \alpha _{i}\right\vert }^{k_{\varepsilon }}$
for every $i\in \{0,1,...,n\}$. Hence $d^{\mu ^{k}}(x_{1},y_{1})\overset{%
\text{(2)}}{\leq }d^{\mu ^{k_{\varepsilon }}}(x_{1},y_{1})\leq \overset{n}{%
\underset{i=0}{\dsum }}z_{\left\vert \alpha _{i}\right\vert }<d^{\mu
}(x_{1},y_{1})+\varepsilon $, so $0\leq d^{\mu ^{k}}(x_{1},y_{1})-d^{\mu
}(x_{1},y_{1})<\varepsilon $ for every $k\in \mathbb{N}$, \thinspace $k\geq
k_{\varepsilon }$, i.e. $d^{\mu }(x_{1},y_{1})=\underset{k\rightarrow \infty 
}{\lim }d^{\mu ^{k}}(x_{1},y_{1})$.

Finally $d^{\mu }(x,y)\overset{\text{(3)}}{\geq }M-\underset{k=0}{\overset{%
\infty }{\dsum }}(d^{\mu ^{k}}(x,y)-d^{\mu ^{k+1}}(x,y))\overset{\text{(1)}}{%
\geq }M-\underset{k=0}{\overset{\infty }{\dsum }}\frac{M}{2^{k+2}}=\frac{M}{2%
}>0$, so $d^{\mu }(x,y)>0$.

iii) We have $d^{\mu }=\underset{k\in \mathbb{N}}{\inf }d^{\mu ^{k}}\leq
d^{\mu 0}=M$. $\square $

\bigskip

\textbf{A bounded and complete metric }$d$\textbf{\ on }$X$

\bigskip

\textbf{Proposition 3.13.} \textit{In the above framework, there exists a
sequence }$(\mu _{n})_{n\in \mathbb{N}}$\textit{\ of decreasing sequences
such that the function }$\rho :X\times X\rightarrow \lbrack 0,\infty )$, 
\textit{given by }$\rho (x,y)=\underset{n=0}{\overset{\infty }{\dsum }}\frac{%
1}{2^{n}}d^{\mu _{n}}(x,y)$\textit{\ for every }$x,y\in X$\textit{, is a
bounded metric.}

\textit{Proof}. Let us consider a fixed $M>0$. For every $x,y\in A$, $x\neq
y $, based on Proposition 3.12, there exists a decreasing sequence $\mu
_{x,y}=(z_{n}^{x,y})_{n\in \mathbb{N}}$ such that $d^{\mu _{x,y}}(x,y)>0$, $%
\underset{n\rightarrow \infty }{\lim }z_{n}^{x,y}=0$ and $d^{\mu _{x,y}}\leq
M$. In the sequel, by $\tau _{A}$ we mean the topology on $A$ that was
defined on the proof of Theorem 3.4, while by $d$ we mean the metric on $A$
given by the same result.

\textit{Claim 1.} \textit{In the above framework, there exist }$%
D_{x},D_{y}\in \tau _{A}$\textit{\ such that:}

\textit{i) }$x\in D_{x}$\textit{\ and }$y\in D_{y}$\textit{;}

\textit{ii) }$d^{\mu _{x,y}}(u,v)>0$\textit{\ for every }$u\in D_{x}$\textit{%
\ and every }$v\in D_{y}$\textit{.}

\textit{Justification of claim 1}. Since $\underset{n\rightarrow \infty }{%
\lim }z_{n}^{x,y}=0$, we can choose $n\in \mathbb{N}$ such that $z_{n}^{x,y}<%
\frac{d^{\mu _{x,y}}(x,y)}{4}$. Then $D_{x}\overset{not}{=}\underset{%
\left\vert \alpha \right\vert =n,x\in A_{\alpha }}{\cup }A_{\alpha }\in \tau
_{A}$ (since, according to the observation made after Claim 1 from the proof
of Theorem 3.4, we have $A=\underset{\left\vert \alpha \right\vert =n}{\cup }%
A_{\alpha }=\underset{\left\vert \alpha \right\vert =n,x\in A_{\alpha }}{%
\cup }A_{\alpha }\cup \underset{\left\vert \alpha \right\vert =n,x\notin
A_{\alpha }}{\cup }A_{\alpha }$ and $\underset{\left\vert \alpha \right\vert
=n,x\notin A_{\alpha }}{\cup }A_{\alpha }$ is compact\ as a finite union of
compact sets) and $D_{y}\overset{not}{=}\underset{\left\vert \alpha
\right\vert =n,y\in A_{\alpha }}{\cup }A_{\alpha }\in \tau _{A}$ (the same
argument).\ For $u\in D_{x}$ and $v\in D_{y}$, since $d^{\mu
_{x,y}}(x,u)\leq z_{n}^{x,y}<\frac{d^{\mu _{x,y}}(x,y)}{4}$ and $d^{\mu
_{x,y}}(y,v)\leq z_{n}^{x,y}<\frac{d^{\mu _{x,y}}(x,y)}{4}$, we have $d^{\mu
_{x,y}}(u,v)\geq d^{\mu _{x,y}}(x,y)-d^{\mu _{x,y}}(y,v)-d^{\mu
_{x,y}}(x,u)\geq d^{\mu _{x,y}}(x,y)-2\frac{d^{\mu _{x,y}}(x,y)}{4}=\frac{%
d^{\mu _{x,y}}(x,y)}{2}>0$ and the justification of the claim in done.

Hence, for every $\varepsilon >0$, from the open cover (provided by Claim
1)\linebreak\ $(D_{x}\times D_{y})_{x,y\in A\times A}$ of the compact set $%
K_{\varepsilon }\overset{not}{=}\{(x,y)\in A\times A\mid d(x,y)\geq
\varepsilon \}$ we can extract a finite open cover, so there exist the
decreasing sequences $\mu ^{1},...,\mu ^{p_{\varepsilon }}$ such that for
every $(x,y)\in K_{\varepsilon }$ there exists $j_{x,y}\in
\{1,2,...,p_{\varepsilon }\}$ having the property that $d^{\mu
^{j_{x,y}}}(x,y)>0$. Consequently, as $\{(x,y)\in A\times A\mid x\neq y\}=%
\underset{n\in \mathbb{N}}{\cup }K_{\frac{1}{n}}$, there exists a sequence $%
(\mu _{n})_{n\in \mathbb{N}}$\ of decreasing sequences such that\textit{\ }%
for every $(x,y)\in A\times A$, $x\neq y$, we can find $n_{x,y}\in \mathbb{N}
$ having the property that $d^{\mu _{n_{x,y}}}(x,y)>0$. Moreover 
\begin{equation}
d^{\mu _{n}}\leq M  \tag{1}
\end{equation}%
and 
\begin{equation}
\underset{k\rightarrow \infty }{\lim }z_{k}^{n}=0\text{,}  \tag{2}
\end{equation}%
for every $n\in \mathbb{N}$, where $\mu _{n}=(z_{k}^{n})_{k\in \mathbb{N}}$.

Now we define the function $\rho :X\times X\rightarrow \lbrack 0,\infty )$
by $\rho (x,y)=\underset{n=0}{\overset{\infty }{\dsum }}\frac{1}{2^{n}}%
d^{\mu _{n}}(x,y)$\ for every $x,y\in X$. As $d^{\mu _{n}}\overset{\text{(1)}%
}{\leq }M$ for every $n\in \mathbb{N}$, $\rho $ is well defined and,
moreover, $\rho \leq 2M$ for every $x,y\in X$, i.e. $\rho $ is bounded. It
is clear that $\rho (x,x)=0$, $\rho (x,y)=\rho (y,x)$ and $\rho (x,y)\leq
\rho (x,z)+\rho (z,y)$ for every $x,y,z\in X$. Moreover, $\rho
(x,y)=0\Rightarrow x=y$ for every $x,y\in X$. Indeed, for $x,y\in X$, $x\neq
y$, we divide the discussion into the following cases: a) $x,y\in A$; b) $%
x\in X\smallsetminus A$ or $y\in X\smallsetminus A$. In case a), we can find 
$n_{x,y}\in \mathbb{N}$ having the property that $d^{\mu _{n_{x,y}}}(x,y)>0$%
, so $\rho (x,y)=\underset{n=0}{\overset{\infty }{\dsum }}\frac{1}{2^{n}}%
d^{\mu _{n}}(x,y)\geq \frac{1}{2^{_{n_{x,y}}}}d^{\mu _{_{n_{x,y}}}}(x,y)>0$,
hence $\rho (x,y)>0$. In situation b), from Proposition 3.8, d) , we infer
that $d^{\mu _{n}}(x,y)>0$ for every $n\in \mathbb{N}$, so $\rho (x,y)>0$.
We conclude that $\rho $ is a metric. $\square $

\bigskip

In the above framework, we consider the sequence $\eta =(z_{k})_{k\in 
\mathbb{N}}$, where $z_{k}=\underset{n=0}{\overset{\infty }{\sum }}\frac{1}{%
2^{n}}z_{k}^{n}$.

\bigskip

\textbf{Proposition 3.14 }(The properties of the sequence $\eta $). \textit{%
In the above framework, the sequence }$\eta $\textit{\ has the following
properties:}

\textit{a) it is well define;}

\textit{b) it is decreasing;}

\textit{c) }$\underset{k\rightarrow \infty }{\lim }z_{k}=0$\textit{.}

\textit{Proof}.

a) As the series $\underset{n=0}{\overset{\infty }{\sum }}\frac{1}{2^{n}}$
is convergent and $z_{k}^{n}\leq M$ for every $k,n\in \mathbb{N}$, the
comparison test yields the conclusion.

b) As $z_{k+1}^{n}\leq z_{k}^{n}$ for every $k,n\in \mathbb{N}$, the same
comparison test assures us that $z_{k+1}\leq z_{k}$ for every $k\in \mathbb{N%
}$.

c) Let us consider an arbitrary $\varepsilon >0$. Since $\underset{%
k\rightarrow \infty }{\lim }z_{k}^{0}=\underset{k\rightarrow \infty }{\lim }%
\frac{1}{2}z_{k}^{1}=...=\underset{k\rightarrow \infty }{\lim }\frac{1}{%
2^{n_{\varepsilon }-1}}z_{k}^{n_{\varepsilon }-1}=0$, where $n_{\varepsilon
}=3+[\log _{2}\frac{M}{\varepsilon }]$, there exists $k_{\varepsilon }\in 
\mathbb{N}$ such that $0<z_{k}^{0}<\frac{\varepsilon }{2n_{\varepsilon }}$, $%
0<\frac{1}{2}z_{k}^{1}<\frac{\varepsilon }{2n_{\varepsilon }}$, ..., $0<%
\frac{1}{2^{n_{\varepsilon }-1}}z_{k}^{n_{\varepsilon }-1}<\frac{\varepsilon 
}{2n_{\varepsilon }}$ for every $k\in \mathbb{N}$, $k\geq k_{\varepsilon }$.
Consequently we have $0\leq z_{k}=z_{k}^{0}+\frac{1}{2}z_{k}^{1}+...+\frac{1%
}{2^{n_{\varepsilon }-1}}z_{k}^{n_{\varepsilon }-1}+\frac{1}{%
2^{n_{\varepsilon }}}z_{k}^{n_{\varepsilon }}+\frac{1}{2^{n_{\varepsilon }+1}%
}z_{k}^{n_{\varepsilon }+1}+...+\frac{1}{2^{n}}z_{k}^{n}+...\leq
n_{\varepsilon }\frac{\varepsilon }{2n_{\varepsilon }}+M(\frac{1}{%
2^{n_{\varepsilon }}}+\frac{1}{2^{n_{\varepsilon }+1}}+...+\frac{1}{2^{n}}%
+...)=\frac{\varepsilon }{2}+\frac{M}{2^{n_{\varepsilon }-1}}<\frac{%
\varepsilon }{2}+\frac{\varepsilon }{2}=\varepsilon $ for every $k\in 
\mathbb{N}$, $k\geq k_{\varepsilon }$ and the conclusion follows. $\square $

\bigskip

Now we can consider the semi-metric $d^{\eta }\overset{not}{=}\delta $.

\bigskip

\textbf{Proposition 3.15 }(The properties of the metric $\delta $). \textit{%
In the above framework, }$\delta $\textit{\ has the following properties:}

\textit{a) }$\rho \leq \delta $\textit{;}

\textit{b) }$\delta \leq 2M$\textit{;}

\textit{c) }$(X,\delta )$\textit{\ is a bounded and complete metric space.}

\textit{Proof}.

a) For $x,y\in X$, $x\neq y$, $p\in \mathbb{N}$ and $\alpha _{0},\alpha
_{1},...,\alpha _{p}\in \Lambda ^{\ast }(I)$ such that $x\in \overset{\sim }{%
X_{\alpha _{0}}}$, $y\in \overset{\sim }{X_{\alpha _{n}}}$ and $\overset{%
\sim }{X_{\alpha _{i}}}\cap \overset{\sim }{X_{\alpha _{i+1}}}\neq \emptyset 
$ for every $i\in \{0,1,...,p-1\}$, we have $\frac{1}{2^{n}}d^{\mu
_{n}}(x,y)\leq \frac{1}{2^{n}}\underset{i=0}{\overset{p}{\sum }}%
z_{\left\vert \alpha _{i}\right\vert }^{n}$ for every $n\in \mathbb{N}$, so $%
\underset{n=0}{\overset{\infty }{\sum }}\frac{1}{2^{n}}d^{\mu _{n}}(x,y)\leq 
\underset{n=0}{\overset{\infty }{\sum }}\frac{1}{2^{n}}z_{\left\vert \alpha
_{0}\right\vert }^{n}+...+\underset{n=0}{\overset{\infty }{\sum }}\frac{1}{%
2^{n}}z_{\left\vert \alpha _{p}\right\vert }^{n}=z_{\left\vert \alpha
_{0}\right\vert }+...+z_{\left\vert \alpha _{p}\right\vert }$. Hence $\rho
(x,y)=\underset{n=0}{\overset{\infty }{\sum }}\frac{1}{2^{n}}d^{\mu
_{n}}(x,y)\leq d^{\eta }(x,y)=\delta (x,y)$. As the last inequality is also
true for $x=y$, the justification of a) is done.

b) As $z_{k}^{n}\leq M$ for every $k,n\in \mathbb{N}$, we deduce that $z_{k}=%
\underset{n=0}{\overset{\infty }{\sum }}\frac{1}{2^{n}}z_{k}^{n}\leq M%
\underset{n=0}{\overset{\infty }{\sum }}\frac{1}{2^{n}}=2M$ for every $k\in 
\mathbb{N}$. Hence $\eta \prec \theta $, where $\theta =(y_{k})_{k\in 
\mathbb{N}}$, $y_{k}=2M$ for every $k\in \mathbb{N}$, and we infer that $%
\delta =d^{\eta }\leq d^{\theta }=2M$.

c) Since $\delta (x,y)=0\overset{\text{a)}}{\Rightarrow }\rho (x,y)=0\overset%
{\text{Proposition 3.13}}{\Rightarrow }x=y$ we conclude that $\delta $ is a
metric on $X$. According to b) it is bounded. In order to prove that $%
(X,\delta )$ is complete, let us consider a Cauchy sequence $(x_{n})_{n\in 
\mathbb{N}}$. By passing to a subsequence, we divide the discussion into the
following two cases (see the proof of Proposition 3.8 for the definition of
the function $m$): a) there exists $N\in \mathbb{N}$ such that $m(x_{n})\leq
N$ for every $n\in \mathbb{N}$; b) $\underset{n\rightarrow \infty }{\lim }%
m(x_{n})=\infty $. In the first case, we have $\delta (x_{n},x_{m})=\inf \{%
\overset{p}{\underset{i=0}{\dsum }}z_{\left\vert \alpha _{i}\right\vert
}\mid $there exist $p\in \mathbb{N}$ and $\alpha _{0},\alpha _{1},...,\alpha
_{p}\in \Lambda ^{\ast }(I)$ such that $x_{n}\in \overset{\sim }{X_{\alpha
_{0}}}$, $x_{m}\in \overset{\sim }{X_{\alpha _{p}}}$ and $\overset{\sim }{%
X_{\alpha _{i}}}\cap \overset{\sim }{X_{\alpha _{i+1}}}\neq \emptyset $ for
every $i\in \{0,1,...,p-1\}\}\overset{(z_{n})_{n\in \mathbb{N}}\text{ is
decreasing}}{\geq }z_{N}$ for every $x_{n}\neq x_{m}$, so, as $(x_{n})_{n\in 
\mathbb{N}}$ is Cauchy, there exists $n_{0}\in \mathbb{N}$ such that $%
x_{n_{0}}=x_{n_{0}+1}=x_{n_{0}+2}=...$ and consequently the sequence $%
(x_{n})_{n\in \mathbb{N}}$ is convergent. In the second case, for each $n\in 
\mathbb{N}$, there exists $\alpha _{n}\in \Lambda ^{\ast }(I)$ such that $%
x_{n}\in \overset{\sim }{X_{\alpha _{n}}}$ and $\left\vert \alpha
_{n}\right\vert =m(x_{n})$. Therefore an argument similar to the one used in
the proof of Proposition 3.11 assures us that one can pick $\alpha \in
\Lambda (I)$ such that $[\alpha ]_{n}=x_{n}\in \overset{\sim }{X_{[\alpha
]_{n}}}$ for every $n\in \mathbb{N}$. Hence $\delta (x_{n},a_{\alpha })\leq
z_{m(x_{n})}$ for every $n\in \mathbb{N}$ which implies that the sequence $%
(x_{n})_{n\in \mathbb{N}}$ is convergent (having the limit $a_{\alpha }$). $%
\square $

\bigskip

Now let us consider a fixed strictly increasing sequence $(c_{n})_{n\in 
\mathbb{N}}$ such that $c_{0}=1$, $(\frac{c_{n}}{c_{n+1}})_{n\in \mathbb{N}}$
is strictly increasing and $c_{n}\leq 2$, $\frac{1}{2}\leq \frac{c_{n}}{%
c_{n+1}}$ for every $n\in \mathbb{N}$ and the function $d:X\times
X\rightarrow \lbrack 0,\infty )$ given by $d(x,y)=\underset{\alpha \in
\Lambda ^{\ast }(I)}{\sup }c_{\left\vert \alpha \right\vert }\delta
(f_{\alpha }(x),f_{\alpha }(y))$ for every $x,y\in X$.

\bigskip

\textbf{Proposition 3.16}. \textit{In the above framework, }$(X,d)$\textit{\
is a bounded and complete metric space.}

\textit{Proof}. It follows form the inequality $\delta \leq d\overset{\text{%
Proposition 3.8, e)}}{\leq }2\delta $ and the fact that, according to
Proposition 3.15, c), $(X,\delta )$ is a bounded and complete metric space. $%
\square $

\bigskip

A word of warning: Even though we use the same notation, namely $d$, for the
metric from Theorem 3.4 and for the one from Proposition 3.16, it is clear
that they are different objects, the first one being a distance on $A$,
while the second one is a metric on $X$.

\bigskip

\textbf{A comparison function }$\varphi $\textbf{\ which makes }$\varphi $%
\textbf{-contractions with respect to }$d$\textbf{\ all the functions of the
family having attractor}

\bigskip

\textbf{Lemma. 3.17.} \textit{In the above framework, we have }$d(X_{\alpha
})\leq d(\overset{\sim }{X_{\alpha }})\leq 2z_{\left\vert \alpha \right\vert
}$\textit{\ for every }$\alpha \in \Lambda ^{\ast }(I)$\textit{.}

\textit{Proof}. For every $x,y\in \overset{\sim }{X_{\alpha }}$ and $\beta
\in \Lambda ^{\ast }(I)$ we have $\delta (f_{\beta }(x),f_{\beta }(y))$%
\linebreak $\overset{f_{\beta }(x),f_{\beta }(y)\overset{\text{Proposition
3.7, e)}}{\in }\overset{\sim }{X_{\beta \alpha }}}{\leq }z_{\left\vert \beta
\alpha \right\vert }\leq z_{\left\vert \alpha \right\vert }$, so $d(x,y)=%
\underset{\beta \in \Lambda ^{\ast }(I)}{\sup }c_{\left\vert \beta
\right\vert }\delta (f_{\beta }(x),f_{\beta }(y))\leq 2z_{\left\vert \alpha
\right\vert }$ and consequently $d(X_{\alpha })\leq d(\overset{\sim }{%
X_{\alpha }})\leq 2z_{\left\vert \alpha \right\vert }$\textit{\ }for every $%
\alpha \in \Lambda ^{\ast }(I)$. $\square $

\bigskip

Let us fix $M>0$. Taking into account Propostion 3.14, c) and Lemma 3.17,
for every $r\in (0,4M]$, there exists $n_{r}\in \mathbb{N}$ such that $%
d(X_{\alpha })\leq d(\overset{\sim }{X_{\alpha }})\leq \frac{r}{20}$ for
every $\alpha \in \Lambda ^{\ast }(I)$ with the property that $\left\vert
\alpha \right\vert \geq n_{r}$. For every $r\in (0,4M)$ we consider the
comparison function $\varphi _{r}:[0,\infty )\rightarrow \lbrack 0,\infty )$%
, given by $\varphi _{r}(x)=\{%
\begin{array}{cc}
0\text{,} & x\in \lbrack 0,r-\rho _{r}) \\ 
\frac{c_{n_{r}}}{c_{n_{r}}+1}x\text{,} & x\in \lbrack r-\rho _{r},r+\rho
_{r}] \\ 
\frac{c_{n_{r}}}{c_{n_{r}}+1}(r+\rho _{r})\text{,} & x\in (r+\rho
_{r},\infty )%
\end{array}%
$, where $\rho _{r}\in (0,\min \{4M-r,\frac{r}{2}\})$. We also consider the
comparison function $\varphi _{4M}:[0,\infty )\rightarrow \lbrack 0,\infty )$%
, given by $\varphi _{4M}(x)=\{%
\begin{array}{cc}
0\text{,} & x\in \lbrack 0,2M) \\ 
\frac{c_{n_{M}}}{c_{nM}+1}x\text{,} & x\in \lbrack 2M,4M] \\ 
\frac{c_{n_{M}}}{c_{n_{M}}+1}4M\text{,} & x\in (4M,\infty )%
\end{array}%
$.

\bigskip

\textbf{Lemma. 3.18}. \textit{In the above framework, we have }$%
d(f_{i}(x),f_{i}(y))\leq \varphi _{r}(d(x,y))$\textit{\ for every }$i\in I$, 
$r\in (0,4M)$ \textit{and }$x,y\in X$\textit{\ having the property that }$%
d(x,y)\in \lbrack r-\rho _{r},r+\rho _{r}]$\textit{. Moreover} $%
d(f_{i}(x),f_{i}(y))\leq \varphi _{4M}(d(x,y))$\textit{\ for every }$i\in I$ 
\textit{and }$x,y\in X$\textit{\ having the property that }$d(x,y)\in
\lbrack 2M,4M]$.

\textit{Proof}. We treat only the situation $r\in (0,4M)$ (the proof for $%
r=4M$ being similar). We divide the discussion into two cases:

a) $d(f_{i}(x),f_{i}(y))<\frac{r}{10}$;

b) $\frac{r}{10}\leq d(f_{i}(x),f_{i}(y))$.

In the first case we have $d(f_{i}(x),f_{i}(y))<\frac{r}{10}<\frac{r}{4}=%
\frac{1}{2}\frac{r}{2}\leq \frac{1}{2}(r-\rho _{r})\leq \frac{1}{2}%
d(x,y)\leq \frac{c_{n_{r}}}{c_{n_{r}+1}}d(x,y)=\varphi _{r}(d(x,y))$.

In the second case, noting that $\delta (f_{\alpha }(f_{i}(x)),f_{\alpha
}(f_{i}(y)))\leq d(f_{\alpha }(f_{i}(x)),f_{\alpha }(f_{i}(y)))$\linebreak $%
\leq d(X_{\alpha })\leq \frac{r}{20}<\frac{r}{10}$ for every $\alpha \in
\Lambda ^{\ast }(I)$ with the property that $\left\vert \alpha \right\vert
\geq n_{r}$, we conclude that $d(f_{i}(x),f_{i}(y))=\underset{\alpha \in
\Lambda ^{\ast }(I)}{\sup }c_{\left\vert \alpha \right\vert }\delta
(f_{\alpha }(f_{i}(x)),f_{\alpha }(f_{i}(y)))=\underset{\alpha \in \Lambda
^{\ast }(I),\left\vert \alpha \right\vert \leq n_{r}}{\max }c_{\left\vert
\alpha \right\vert }\delta (f_{\alpha }(f_{i}(x)),f_{\alpha }(f_{i}(y)))$,
so there exists $\alpha _{0}\in \Lambda ^{\ast }(I)$ such that $\left\vert
\alpha _{0}\right\vert \leq n_{r}$ and 
\begin{equation}
d(f_{i}(x),f_{i}(y))=c_{\left\vert \alpha _{0}\right\vert }\delta (f_{\alpha
_{0}}(f_{i}(x)),f_{\alpha _{0}}(f_{i}(y)))\text{.}  \tag{1}
\end{equation}%
As $c_{\left\vert \alpha _{0i}\right\vert }\delta (f_{\alpha
_{0i}}(x),f_{\alpha _{0}i}(y))\leq d(x,y)$, using $(1)$, we get $\frac{%
c_{\left\vert \alpha _{0i}\right\vert }}{c_{\left\vert \alpha
_{0}\right\vert }}d(f_{i}(x),f_{i}(y))\leq d(x,y)$, i.e. $%
d(f_{i}(x),f_{i}(y))\leq \frac{c_{\left\vert \alpha _{0}\right\vert }}{%
c_{\left\vert \alpha _{0}\right\vert +1}}d(x,y)\leq \frac{c_{n_{r}}}{%
c_{n_{r}+1}}d(x,y)=\varphi _{r}(d(x,y))$. $\square $

\bigskip

Note that the family consisting of the intervals $(2M,5M)$ and $(r-\rho
_{r},r+\rho _{r})$, where $r\in (0,4M)$, is an open cover of $(0,4M]$ which
is Lindel\"{o}f and paracompact, so there exists a sequence $(r_{n})_{n\in 
\mathbb{N}}$ of elements from $(0,4M]$ such that $(0,4M]=\underset{n\in 
\mathbb{N}}{\cup }[r_{n}-\rho _{r_{n}},r_{n}+\rho _{r_{n}}]$ and the family $%
\{[r_{n}-\rho _{r_{n}},r_{n}+\rho _{r_{n}}]\}_{n\in \mathbb{N}}$ is locally
finite, where by $[r_{n}-\rho _{r_{n}},r_{n}+\rho _{r_{n}}]$ we mean $%
[2M,4M] $ in case that $r_{n}=4M$.

\bigskip

\textbf{Lemma 3.19}. \textit{In the above framework, the function }$\varphi =%
\underset{n\in \mathbb{N}}{\sup }\varphi _{r_{n}}$\textit{\ is a comparison
function.}

\textit{Proof}. It is obvious that $\varphi $ is increasing and that $%
\varphi (t)<t$ for every $t>0$ since all the functions $\varphi _{r_{n}}$
have these properties. Moreover, for every $t\in \lbrack 0,\infty )$ there
exists a neighborhood $V_{t}$ of $t$ which intersects only a finite number
of intervals $[r_{n}-\rho _{r_{n}},r_{n}+\rho _{r_{n}}]$ and consequently $%
\varphi _{\mid V_{t}}$ is continuous since it can be presented as the
maximum of a finite set of continuous functions. Hence $\varphi $ is
continuous. $\square $

\bigskip

\textbf{Lemma 3.20}. \textit{In the above framework, all the functions }$%
f_{i}$\textit{\ are }$\varphi $\textit{-contractions.}

\textit{Proof}. For $x,y\in X$, $x\neq y$, we have $d(x,y)\in (0,4M]=%
\underset{n\in \mathbb{N}}{\cup }[r_{n}-\rho _{r_{n}},r_{n}+\rho _{r_{n}}]$,
so there exists $n_{0}\in \mathbb{N}$ such that $d(x,y)\in \lbrack
r_{n_{0}}-\rho _{r_{n_{0}}},r_{n_{0}}+\rho _{r_{n_{0}}}]$ and therefore we
have $d(f_{i}(x),f_{i}(y))\overset{\text{Lemma 3.18}}{\leq }\varphi
_{r_{n_{0}}}(d(x,y))\leq \varphi (d(x,y))$ for every\textit{\ }$i\in I$. As
the last inequality is obviously true for $x=y$, we conclude that $f_{i}$%
\textit{\ }is $\varphi $-contraction for every $i\in I$. $\square $

\bigskip

We summarize the above facts in the following:

\bigskip

\textbf{Theorem 3.21.} \textit{Given a family of functions\ }$(f_{i})_{i\in
I}$\textit{\ having attractor, there exists a metric }$d$\textit{\ on }$X$%
\textit{\ and a comparison function }$\varphi $\textit{\ such that:}

\textit{\qquad a) the metric space }$(X,d)$\textit{\ is complete and bounded;%
}

\textit{\qquad b) }$f_{i}$\textit{\ is }$\varphi $\textit{-contraction with
respect to }$d$\textit{\ for every }$i\in I$\textit{.}

\bigskip

Combining Proposition 3.1 and Theorem 3.21 we obtain the following:

\bigskip

\textbf{Theorem 3.22.} \textit{Given} $(f_{i})_{i\in I}$\textit{\ a family
of functions, where }$f_{i}:X\rightarrow X$\textit{\ and }$I$\textit{\ is
finite, the following two statements are equivalent:}

\textit{I. There exists a metric }$d$\textit{\ on }$X$\textit{\ and a
comparison function }$\varphi $\textit{\ such that:}

\textit{\qquad a) the metric space }$(X,d)$\textit{\ is complete and bounded;%
}

\textit{\qquad b) }$f_{i}$\textit{\ is }$\varphi $\textit{-contraction with
respect to }$d$\textit{\ for every }$i\in I$\textit{.}

\textit{II. The following two statements are valid:}

\textit{\qquad a) For every }$\alpha \in \Lambda (I)$\textit{, the set }$%
\underset{n\in \mathbb{N}}{\cap }X_{[\alpha ]_{n}}$\textit{\ has a unique
element which is denoted by }$a_{\alpha }$\textit{.}

\textit{\qquad b) If }$a_{\alpha }\neq a_{\beta }$\textit{, where }$\alpha
,\beta \in \Lambda (I)$\textit{, then there exists }$n_{0}\in \mathbb{N}$%
\textit{\ such that }$X_{[\alpha ]_{n_{0}}}\cap X_{[\beta
]_{n_{0}}}=\emptyset $\textit{.}

\bigskip

\textbf{4.} \textbf{FINAL\ REMARKS}

\bigskip

\textbf{The unbounded case}

\bigskip

The following result removes the boundedness restriction on the metric $d$.

\bigskip

\textbf{Theorem 4.1.} \textit{Given} $(f_{i})_{i\in I}$\textit{\ a family of
functions, where }$f_{i}:X\rightarrow X$\textit{\ and }$I$\textit{\ is
finite, the following two statements are equivalent:}

\textit{I. There exists a metric }$D$\textit{\ on }$X$\textit{\ and a
comparison function }$\varphi $\textit{\ such that:}

\textit{\qquad a) the metric space }$(X,D)$\textit{\ is complete;}

\textit{\qquad b) }$f_{i}$\textit{\ is }$\varphi $\textit{-contraction with
respect to }$D$\textit{\ for every }$i\in I$\textit{.}

\textit{II. There exists a subset }$X_{1}$\textit{\ of }$X$\textit{\ such
that the following four statements are valid:}

\qquad \textit{a)} $F(X_{1})\subseteq X_{1}$\textit{.}

\textit{\qquad b) For every }$\alpha \in \Lambda (I)$\textit{, the set }$%
\underset{n\in \mathbb{N}}{\cap }(X_{1})_{[\alpha ]_{n}}$\textit{\ has a
unique element which is denoted by }$a_{\alpha }$\textit{.}

\textit{\qquad c) If }$a_{\alpha }\neq a_{\beta }$\textit{, where }$\alpha
,\beta \in \Lambda (I)$\textit{, then there exists }$n_{0}\in \mathbb{N}$%
\textit{\ such that }$(X_{1})_{[\alpha ]_{n_{0}}}\cap (X_{1})_{[\beta
]_{n_{0}}}=\emptyset $\textit{.}

\qquad \textit{d) For every }$x\in X$\textit{\ there exists} $n_{x}\in 
\mathbb{N}$ \textit{such that }$F^{[n_{x}]}(\{x\})\subseteq X_{1}$\textit{,
where }$F:\mathcal{P}(X)\rightarrow \mathcal{P}(X)$\textit{\ is given by }$%
F(C)=\underset{i\in I}{\cup }f_{i}(C)$\textit{\ for every subset }$C$\textit{%
\ of }$X$\textit{.}

\textit{Proof}.

I)$\Rightarrow $II) We choose $X_{1}=B(A,r)$, where $r>0$ and $A$ is the
unique fixed point of the function $F_{\mathcal{S}}$ defined on Remark 3.2,
ii). For the verification of a) we choose $v\in F(B(A,r))$. Then there
exists $i\in I$ and $x\in B(A,r)$ such that $v=f_{i}(x)$ and there exists $%
y\in A$ such that $d(x,y)<r$. Then $d(v,f_{i}(y))=d(f_{i}(x),f_{i}(y))\leq
\varphi (d(x,y))<\varphi (r)$ and consequently, as $f_{i}(y)\in
f_{i}(A)\subseteq A$, we conclude that $v\in B(A,r)$. Hence $%
F(B(A,r))\subseteq B(A,r)$.The properties b) and c) can be proved with
exactly the same techniques used in the proof of Proposition 3.1. The
property d) results from the fact that $\underset{n\rightarrow \infty }{\lim 
}h(F^{[n]}(\{x\}),A)=0$ (see Remark 3.2, ii)).

II)$\Rightarrow $I) According to Theorem 3.21, taking into account b) and
c), there exists a metric $d$\ on $X_{1}$\ and a comparison function $%
\varphi $\ such that: a) the metric space $(X_{1},d)$\ is complete and
bounded; b) $f_{i}$\ is $\varphi $-contraction with respect to $d$\ for
every $i\in I$.

For a given $a\in (0,1)$, the function $\psi :[0,\infty )\rightarrow \lbrack
0,\infty )$ given by 
\begin{equation*}
\psi (t)=\underset{x\in \lbrack 0,t]}{\sup }\{ax+\varphi (t-x)\}=\sup
\{\varphi _{1}(t_{1})+\varphi (t_{2})\mid t_{1},t_{2}\geq 0,t_{1}+t_{2}\leq
t\}\text{,}
\end{equation*}%
where $\varphi _{1}(t)=at$, for every $t\geq 0$, is a comparison function
(see Fact 10 from the proof of Theorem 3.1 from [20]).

Note that, taking into account Remark 2.2, ii), we have $\varphi \leq \psi $
and $\varphi _{1}\leq \psi $.

We consider the function $D:X\times X\rightarrow \lbrack 0,\infty )$ given
by 
\begin{equation*}
D(x,y)=\{%
\begin{array}{cc}
d(x,y)\text{,} & x,y\in X_{1} \\ 
\max \{Ma^{-l(x)},Ma^{-l(y)}\}\text{,} & \{x,y\}\cap (X\smallsetminus
X_{1})\neq \emptyset \text{ and }x\neq y \\ 
0\text{,} & x=y\in X\smallsetminus X_{1}%
\end{array}%
\text{,}
\end{equation*}%
where $l(x)=\{%
\begin{array}{cc}
-\infty \text{,} & x\in X_{1} \\ 
\min \{n\in \mathbb{N}\mid F^{[n]}(\{x\})\subseteq X_{1}\}\text{,} & x\in
X\smallsetminus X_{1}%
\end{array}%
$ and $M$ is an upper bound for $d$. Note that, according to d), $l(x)\in 
\mathbb{N}$ for every $x\in X\smallsetminus X_{1}$. We use the convention
that $a^{\infty }=0$, so $Ma^{-l(x)}=0$ for $x\in X_{1}$. One can routinely
check that $D$ is a metric on $X$.

Moreover, 
\begin{equation*}
D(f_{i}(x),f_{i}(y))\leq \psi (d(x,y))\text{,}
\end{equation*}
for every $i\in I$, $x,y\in X$.

Indeed, if $x,y\in X_{1}$, then $f_{i}(x),f_{i}(y)\overset{a)}{\in }X_{1}$,
so $D(f_{i}(x),f_{i}(y))=d(f_{i}(x),f_{i}(y))\leq \varphi (d(x,y))\leq \psi
(d(x,y))=\psi (D(x,y))$. If $\{x,y\}\cap (X\smallsetminus X_{1})\neq
\emptyset $ and $x\neq y$, then we divide the discussion into three cases:
1. $l(x)=l(y)=1$. 2. $l(x)\geq l(y)>1$. 3. $l(y)\geq l(x)>1$. In the first
case, as $f_{i}(x)\in F^{[l(x)]}(\{x\})\subseteq X_{1},f_{i}(y)\in
F^{[l(y)]}(\{y\})\subseteq X_{1}$, we have $%
D(f_{i}(x),f_{i}(y))=d(f_{i}(x),f_{i}(y))\leq aa^{-1}M=aD(x,y)=\varphi
_{1}(D(x,y))\leq \psi (D(x,y))$. In the second case, note that $%
l(f_{i}(x))=l(x)-1$ and $l(f_{i}(y))=l(y)-1$, so $D(f_{i}(x),f_{i}(y))=\max
\{Ma^{-l(f_{i}(x))},Ma^{-l(f_{i}(y))}\}=aD(x,y)=\varphi _{1}(D(x,y))\leq
\psi (D(x,y))$. The third case is similar with the second one. If $x=y\in
X\smallsetminus X_{1}$ the conclusion is clear. $\square $

\bigskip

\textbf{Some facts about the topological structure of }$(X,d^{\mu })$

\bigskip

In the framework of the third section, let us suppose that $\sigma $ is a
distance on $X$ such that there exist $(c_{n})_{n\in \mathbb{N}}$ and $%
(d_{n})_{n\in \mathbb{N}}$ having the following properties:

a) $\underset{n\rightarrow \infty }{\lim }c_{n}=\underset{n\rightarrow
\infty }{\lim }d_{n}=0$;

b) $\sigma (x,y)\geq c_{m(x)}$ for every $x,y\in X$, $x\neq y$, with the
convention that $c_{\infty }=0$ (for the definition of $m(x)$ see the proof
of Theorem 3.8, d));

c) $d(X_{\alpha })\leq d_{\left\vert \alpha \right\vert }$ for every $\alpha
\in \Lambda ^{\ast }(I)$.

Let us denote by $\tau $ the topology induced by $\sigma $.

\bigskip

Then one can easily check the following properties:

i) the sets $\overset{\sim }{X_{\alpha }}$ are closed with respect to $\tau $%
;

ii) $\{x\}$ is open with respect to $\tau $ for every $x\in X\smallsetminus
A $;

iii) $(V_{x,n})_{n\in \mathbb{N}^{\ast }}$ is a neighborhood basis for $x$
with respect to $\tau $, where $V_{x,n}=\underset{\alpha \in \Lambda ^{\ast
}(I),\left\vert \alpha \right\vert =n,x\in X_{\alpha }}{\cup }\overset{\sim }%
{X_{\alpha }}$ for every $x\in A$;

iv) the function $\pi :\Lambda (I)\rightarrow A$, given by $\pi (\alpha
)=a_{\alpha }$ for every $\alpha \in \Lambda (I)$, is continuous with
respect to $\tau $;

v) If $(x_{n})_{n\in \mathbb{N}}$ is a sequence of elements from $X$ and $%
x\in X$, then:

\qquad j) for $x\in X\smallsetminus A$: $\underset{n\rightarrow \infty }{%
\lim }x_{n}=x$ with respect to $\tau $ if and only if there exists $n_{0}\in 
\mathbb{N}$ such that $x_{n}=x$ for every $n\in \mathbb{N}$, $n\geq n_{0}$;

\qquad jj) for $x\in A$: $\underset{n\rightarrow \infty }{\lim }x_{n}=x$
with respect to $\tau $ if and only if for every $m\in \mathbb{N}$ there
exists $n_{m}\in \mathbb{N}$ having the property that for every $n\in 
\mathbb{N}$, $n\geq n_{m}$ there exists $\alpha ^{n}\in \Lambda (I)$ such
that $x=a_{\alpha ^{n}}$ and $x_{n}\in X_{[\alpha ^{n}]_{m}}$ for every $%
n\in \mathbb{N}$;

vi) $(X,\sigma )$ is complete.

\bigskip

Note that if $d^{\mu }$ is a distance, where $\mu =(\alpha ^{n})_{n\in 
\mathbb{N}}$ for some $\alpha \in (0,1)$, satisfies the requirements imposed
on the metric $\sigma $ from the previous paragraph. Indeed, take $%
c_{n}=d_{n}=\alpha ^{n}$ for every $n\in \mathbb{N}$ and note that a) is
obvious, b) results from the proof of Theorem 3.8 and c) could be obtained
directly from the definition of $d^{\mu }$. Consequently, according to vi), $%
(X,d^{\mu })$ is complete.

\bigskip

\textbf{The particular case of a family consisting of one function}

\bigskip

For the particular case of a family $(f_{i})_{i\in I}$ having the property
that the set $I$ has one element, we obtain the following converse of
Browder's theorem:

\bigskip

\textbf{Proposition 4.2.} \textit{Given a set }$X$\textit{\ and a function }$%
f:X\rightarrow X$\textit{\ such that }$\underset{n\in \mathbb{N}}{\cap }%
f^{[n]}(X)$\textit{\ is a singleton, there exist a bounded and complete
metric }$d$\textit{\ on }$X$\textit{\ and a comparison function }$\varphi $%
\textit{\ such that }$d(f(x),f(y))\leq \varphi (d(x,y))$\textit{\ for every }%
$x,y\in X$.

\textit{Proof}. $\mathcal{F}=\{f\}$ is a family of functions having
attractor since the second condition from the definition of such a system is
obviously valid as $\Lambda (I)$ has just one element and therefore the
attractor has just one element. Then just apply Theorem 3.21. $\square $

\bigskip

Moreover, the following stronger results is valid (see Theorem 5 from [10]):

\bigskip

\textbf{Proposition 4.3}. \textit{Given a set }$X$\textit{, }$\alpha \in
(0,1)$\textit{\ and a function }$f:X\rightarrow X$\textit{\ such that}$%
\underset{n\in \mathbb{N}}{\cap }f^{[n]}(X)$\textit{\ is a singleton, there
exists a complete and bounded metric }$d$\textit{\ on }$X$\textit{\ such
that }$d(f(x),f(y))\leq \alpha d(x,y)$\textit{\ for every }$x,y\in X$\textit{%
.}

\textit{Proof}. $\mathcal{F}=\{f\}$ is a family of functions having
attractor consisting of just one element. Hence, given $\alpha \in (0,1)$,
Proposition 3.8, f) and g), assures us that $d^{\mu }$ is a bounded distance
and the same line of arguments used in the proof of Proposition 3.8, e),
confirms that $d^{\mu }(f(x),f(y))\leq \alpha d^{\mu }(x,y)$\ for every $%
x,y\in X$, where $\mu =(\alpha ^{n})_{n\in \mathbb{N}}$. Moreover, according
to the note from the end of the previous section, $(X,d^{\mu })$ is
complete. Now just take $\delta =d^{\mu }$. $\square $

\bigskip

Note that \textit{the condition that }$\underset{n\in \mathbb{N}}{\cap }%
f^{[n]}(X)$\textit{\ is a singleton (i.e. there exists a unique }$x_{0}\in X$%
\textit{\ such that }$\underset{n\in \mathbb{N}}{\cap }f^{[n]}(X)=\{x_{0}\}$%
\textit{) implies that }$x_{0}$\textit{\ is the unique fixed point of }$%
f^{[k]}$\textit{\ for every }$k\in \mathbb{N}$.

Indeed, if $\underset{n\in \mathbb{N}}{\cap }f^{[n]}(X)=\{x_{0}\}$, then $%
f^{[k]}(x_{0})\in \underset{n\in \mathbb{N}}{\cap }f^{[n]}(X)$, so $%
f^{[k]}(x_{0})=x_{0}$, i.e. $x_{0}$ is a fixed point of $f^{[k]}$. Moreover,
if $x_{1}\in X$ is a fixed point of $f^{[k]}$, then $x_{1}\in \underset{n\in 
\mathbb{N}}{\cap }f^{[n]}(X)=\{x_{0}\}$, so $x_{1}=x_{0}$ and consequently $%
x_{0}$ is the unique fixed point of $f^{[k]}$.

\bigskip

\textbf{The particular case of a family of functions having attractor with
common fixed point}

\bigskip

Note that \textit{each of the functions of a family of functions\ having
attractor has a unique fixed point}.

Indeed, let us consider $(f_{i})_{i\in I}$ a\textit{\ }family of functions\
having attractor. Then, according to the property a) from the definition of a%
\textit{\ }family of functions\ having attractor, $\underset{n\in \mathbb{N}}%
{\cap }f_{i}^{[n]}(X)=\underset{n\in \mathbb{N}}{\cap }X_{[\theta ]_{n}}$ is
a singleton and if $\underset{n\in \mathbb{N}}{\cap }f_{i}^{[n]}(X)=\{x_{i}%
\} $, then $x_{i}$ is the unique fixed point of $f_{i}$ for every $i\in I$.
Here $\theta $ is the element of $\Lambda (I)$ having all letters equal to $%
i $.

\bigskip

The following proposition is a companion of the result due to Wong (see
[24]) that extends Bessaga's theorem for a finite family of commuting
functions with common unique fixed point. Note that the commutativity of the
family's functions is not part of the hypotheses of our result.

\bigskip

\textbf{Proposition 4.4}. \textit{Given a set }$X$\textit{, }$\alpha \in
(0,1)$\textit{\ and a family of functions\ }$(f_{i})_{i\in I}$\ \textit{%
having attractor, there exists a complete and bounded metric }$d$\textit{\
on }$X$\textit{\ such that }$d(f_{i}(x),f_{i}(y))\leq \alpha d(x,y)$\textit{%
\ for every }$x,y\in X$ \textit{and every }$i\in I$\textit{, provided that
there exists }$x_{0}\in X$\textit{\ such that }$f_{i}(x_{0})=x_{0}$\textit{.}

\textit{Proof}. We have $x_{0}=f_{[\beta ]_{n}}(x_{0})=f_{[\gamma
]_{n}}(x_{0})\in X_{[\beta ]_{n}}\cap X_{[\gamma ]_{n}}$, so $X_{[\beta
]_{n}}\cap X_{[\gamma ]_{n}}\neq \emptyset $ for every $n\in \mathbb{N}$ and
every $\beta ,\gamma \in \Lambda (I)$. Based on the conditions from the
definition of a family of functions\ having attractor, we infer that the
attractor of $(f_{i})_{i\in I}$ has just one element and the same arguments
used in the proof of Proposition 4.3 assure us that for the complete and
bounded metric $d=d^{\mu }$, where $\mu =(\alpha ^{n})_{n\in \mathbb{N}}$,
we have $d(f_{i}(x),f_{i}(y))\leq \alpha d(x,y)$\textit{\ }for every $x,y\in
X$ and every $i\in I$. $\square $

\bigskip

\textbf{References}

\bigskip

[1] M. Barnsley, Fractals Everywhere, Academic Press, Boston, MA, 1988.

[2] V. Berinde, Iterative approximation of fixed points, 2nd Ed., Springer
Verlag, Berlin-Heidelberg-New York, 2007.

[3] C. Bessaga, On the converse of the Banach fixed point principle, Colloq.
Math., \textbf{7} (1959), 41-43.

[4] D.W. Boyd and J.S. Wong, On nonlinear contractions, Proc. Amer Math.
Soc., \textbf{20} (1969), 458-464.

[5] F.E. Browder, On the convergence of successive approximations for
nonlinear functional equations, Indag. Math., \textbf{30} (1968), 27--35.

[6] J. Dugundji, Topology, Allyn and Bacon, Boston, 1968.

[7] D. Dumitru, Attractors of infinite iterated function systems containing
contraction type functions, An. \c{S}tiin\c{t}. Univ. Al. I. Cuza Ia\c{s}i,
Ser. Nou\u{a}, Mat., \textbf{59} (2013), 281-298.

[8] G. Gw\'{o}\'{z}d\'{z}-\L ukawska and J. Jachymski, The
Hutchinson-Barnsley theory for infinite iterated function systems, Bull.
Aust. Math. Soc., \textbf{72} (2005), 441-454.

[9] J.E. Hutchinson, Fractals and self similarity, Indiana Univ. Math. J., 
\textbf{30\ }(1981), 713-747.

[10] J. R. Jachymski, A short proof of the converse to the contraction
principle and some related results, Topol. Methods Nonlinear Anal., \textbf{%
15} (2000),179--186.

[11] J. R. Jachymski, Around Browder's fixed point theorem for contractions,
J. Fixed Point Theory Appl., \textbf{5} (2009), 47-61.

[12] L. Jano\v{s}, A converse of the Banach's contraction theorem, Proc.
Amer. Math. Soc., \textbf{18} (1967), 287--289.

[13] L. Jano\v{s}, An application of combinatorial techniques to a
topological problem, Bull. Austral. Math. Soc., \textbf{9} (1973), 439--443.

[14] A. Kameyama, Distances on topological self-similar sets and the
kneading determinants, J. Math. Kyoto Univ., \textbf{40} (2000), 601-672.

[15] S. Leader, A topological characterization of Banach contractions, Pac.
Jour. Math., \textbf{69} (1977), 461--466.

[16] J. Matkowski, Fixed point theorems for mappings with a contractive
iterate at a point, Proc. Amer. Math. Soc., \textbf{62} (1977), 344-348.

[17] J. Matkowski, Integrable solutions of functional equations,
Dissertations Math. (Rozprawy), 127 (1976).

[18] P. R. Meyers, A converse to Banach's contraction theorem, J. Research
Nat. Bureau of Standards - B. Math. and Math. Physics, \textbf{71B} (1967),
73-76.

[19] R. Miculescu and A. Mihail, On a question of A. Kameyama concerning
self-similar metrics, J. Math. Anal. Appl., \textbf{422 }(2015), 265-271.

[20] R. Miculescu and A. Mihail, A sufficient condition for a finite family
of continuous functions to be transformed into $\psi $-contractions, Ann.
Acad. Sci. Fenn., Math., \textbf{41} (2016), 51-65.

[21] I. A. Rus, Generalized $\varphi $-contractions, Math., Rev. Anal. Num%
\'{e}r. Th\'{e}or. Approximation, Math., \textbf{47} (1982), 175-178.

[22] F. Strobin and J. Swaczyna, On a certain generalisation of the iterated
function system, Bull. Aust. Math. Soc., \textbf{87} (2013), 37-54.

[23] F. Strobin, An application of a fixed point theorem for multifunctions
in a problem of connectedness of attractors of families of IFSs, Bull. Soc.
Sci. Lett. \L \'{o}d\'{z}, S\'{e}r. Rech. D\'{e}form.,\textbf{\ 64 }(2014),
81-93.

[24] J. S. W. Wong, Generalizations of the converse of the contraction
mapping principle, Canad. J. Math., \textbf{18} (1966), 1095-1104.

\end{document}